\documentclass[pdflatex,sn-mathphys-ay]{sn-jnl}
\usepackage{graphicx}%
\usepackage{multirow}%
\usepackage{amsmath,amssymb,amsfonts}%
\usepackage{amsthm}%
\usepackage{mathrsfs}%
\usepackage[title]{appendix}%
\usepackage{xcolor}%
\usepackage{textcomp}%
\usepackage{manyfoot}%
\usepackage{booktabs}%
\usepackage{algorithm}%
\usepackage{algorithmicx}%
\usepackage{algpseudocode}%
\usepackage{listings}%
\usepackage{doi}
\usepackage{breqn}
\usepackage{footnote}
\usepackage{tabularray}
\usepackage{tikz}
\usetikzlibrary{calc, patterns}
\usepackage{multicol}
\usepackage{subcaption}
\usepackage{wrapfig}
\usepackage{tabularx}
\usepackage{threeparttable}
\usepackage{natbib}

\theoremstyle{thmstyleone}%
\theoremstyle{thmstyletwo}%

\theoremstyle{thmstylethree}%
\newtheorem{definition}{Definition}%
\begin{document}

\title[Article Title]{Open Loop Layout Optimization: Feasible Path Planning and Exact Door-to-Door Distance Calculation}

\author*[1]{\fnm{Seyed Mahdi} \sur{Shavarani}}\email{m.shavarani@kent.ac.uk}

\author[2]{\fnm{B\'ela} \sur{Vizv\'ari}}\email{bela.vizvari@emu.edu.tr}

\author[3]{\fnm{Kov\'acs} \sur{Gergely}}\email{kovacs.gergely@edutus.hu}

\affil*[1]{\orgdiv{Centre for Logistics and Sustainability Analytics}, \orgname{Kent Business School, University of Kent}, \country{UK}}

\affil[2]{\orgdiv{Department of Industrial Engineering}, \orgname{Eastern Mediterranean University}, \country{Turkey}}

\affil[3]{\orgdiv{Edutus University}, \country{Hungary}}

\abstract{The Open Loop Layout Problem (OLLP) seeks to position rectangular cells of varying dimensions on a plane without overlap, minimizing transportation costs computed as the flow-weighted sum of pairwise distances between cells. A key challenge in OLLP is to compute accurate inter-cell distances along feasible paths that avoid rectangle intersections. Existing approaches approximate inter-cell distances using centroids, a simplification that can ignore physical constraints, resulting in infeasible layouts or underestimated distances.
This study proposes the first mathematical model that incorporates exact door-to-door distances and feasible paths under the Euclidean metric, with cell doors acting as pickup and delivery points. Feasible paths between doors must either follow rectangle edges as corridors or take direct, unobstructed routes. To address the NP-hardness of the problem, we present a metaheuristic framework with a novel encoding scheme that embeds exact path calculations. Experiments on standard benchmark instances confirm that our approach consistently outperforms existing methods, delivering superior solution quality and practical applicability. }


\maketitle

\section{Introduction}
Layout optimization is a fundamental problem in operational research, focusing on the arrangement of spatial entities within a bounded region to optimize specific objectives, such as minimizing transportation costs or maximizing space utilization, while adhering to constraints like non-overlapping placement and spatial relationships \citep{wan2022differential}. The spatial entities may represent machines in manufacturing, locations in logistics, infrastructure in urban planning, chips in semiconductor design, and assets in facility management, where optimized layouts enhance productivity, streamline resource utilization, and improve process flows \citep{ pouPieMar2021integrating}.

\begin{figure}[ht]
    \centering
    \begin{tikzpicture}[remember picture, overlay]
        \begin{scope}[shift={(-6.4,0)}] 
    \draw[thick] (0,0) rectangle ++(0.5,0.25);
    \node[anchor=west, font=\scriptsize] at (0.6,0.125) {Rectangular Cell};

    \draw[thick] (0,-.4) -- ++(0.5,0);
    \node[anchor=west, font=\scriptsize] at (0.6,-.4) {Transportation Track};
        \end{scope}
    \end{tikzpicture}

    \begin{subfigure}[t]{0.23\textwidth}
        \centering
        \begin{tikzpicture}[scale=0.4]
            \def\rectwidth{0.325}
            \def\rectheight{0.8125}
            \draw[thick] (0,1.5) -- (4.5,1.5);
            \foreach \i in {0, 1.5, 3} {
                \draw[thick] (\i,1.5) rectangle ++(\rectwidth,\rectheight);
            }
        \end{tikzpicture}
        \caption{Spine}
    \end{subfigure}
    \begin{subfigure}[t]{0.23\textwidth}
        \centering
        \begin{tikzpicture}[scale=0.4]
            \def\rectwidth{0.325}
            \def\rectheight{0.8125}
            \draw[thick] (0,1.5) rectangle (4.5,3.39);
            \foreach \x/\y in {0/1.5, 3.5/1.5, 3.5/2.5, 0/2.5} {
                \draw[thick] (\x,\y) rectangle ++(\rectwidth,\rectheight);
            }
        \end{tikzpicture}
        \caption{Closed Loop}
    \end{subfigure}
    \begin{subfigure}[t]{0.23\textwidth}
        \centering
        \begin{tikzpicture}[scale=0.4]
            \def\rectwidth{0.325}
            \def\rectheight{0.8125}
            \foreach \x in {0, 4.5} {
                \draw[thick] (\x,1) -- ++(0,4);
            }
            \foreach \y in {1,2.5,4} {
                \draw[thick] (0,\y) -- (4.5,\y);
            }
            \foreach \x/\y in {0/1, 1.5/1, 3/2.5, 1.5/2.5} {
                \draw[thick] (\x,\y) rectangle ++(\rectwidth,\rectheight);
            }
        \end{tikzpicture}
        \caption{Ladder}
    \end{subfigure}
    \begin{subfigure}[t]{0.23\textwidth}
        \centering
        \begin{tikzpicture}[scale=0.4]
            \def\rectwidth{0.325}
            \def\rectheight{0.8125}
            \foreach \x/\y in {0.5/1.5, 2.5/2, 1.5/3.5, 3.5/2.5} {
                \draw[thick] (\x,\y) rectangle ++(\rectwidth,\rectheight);
            }
        \end{tikzpicture}
        \caption{Open Loop}
    \end{subfigure}
    \caption{Examples of different layout configurations.}
    \label{fig:layout_samples}
\end{figure}
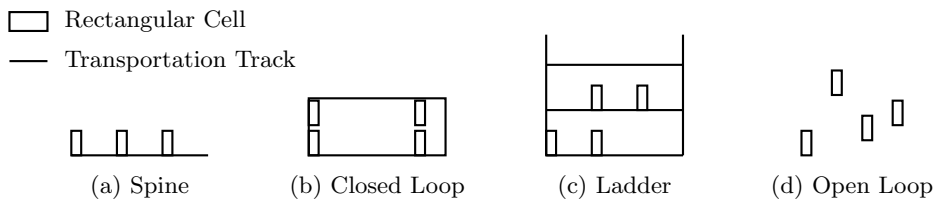

Layout problems are typically categorised into four main types: \textbf{Spine}, \textbf{Closed Loop}, \textbf{Ladder}, and \textbf{Open Loop} \citep{drira2006survey, kang2018closed}, with examples of each illustrated in Figure~\ref{fig:layout_samples}. Each layout type addresses specific constraints and operational needs. Spine layouts streamline movement paths along a single linear track, while Closed Loop layouts organize flow within a rectangular loop, often used in manufacturing. Ladder layouts divide spaces into grid-based zones for efficient access, whereas Open Loop layouts - the most flexible - accommodate diverse configurations and unrestricted flow paths without predefined constraints \citep{ Yang2005LayoutDF}.\\

In contrast to Spine and Closed Loop layouts, where distances can be calculated precisely, Open Loop layouts—even those featuring a single loop—incorporate irregular pathways, branching tracks, and intersecting zones that complicate the computation of shortest paths \citep{de1997complexity, niroomand2013mixed}. Incorporating the concept of "doors" enhances the realism of the model but also significantly increases its complexity \citep{carlo2011analysis}. In realistic scenarios, doors represent designated entry and exit points for inter-cell flow within rectangular entities, restricting movement to specific pathways \citep{GU2010539}. This realistic yet complex feature, combined with the need for exact distance calculations, signals the limitations of traditional approximation and centroid-based methods, which often oversimplify spatial relationships and lead to sub-optimal solutions.  

An inferior solution disrupts flow and reduces operational efficiency ~\citep{niroomand2013mixed, drira2006survey}. In Open Loop Layout Problems (OLLP), these inefficiencies arise from the challenge of accurately determining the shortest feasible paths between objects \citep{ niroomand2013mixed}. 
 We address this challenge by  introducing a mathematical model that employs exact feasible path calculations based on the $l_2$ metric. We further introduce a scalable solution method that efficiently solves benchmark problems in the field, demonstrating significant improvements over existing results.  The model incorporates the "door-to-door" concept to enhance realism and employs advanced metaheuristic algorithms for efficient problem-solving. Through extensive experimentation on benchmark problems, we demonstrate significant improvements over existing approaches, providing a robust framework for tackling real-world layout optimization challenges.

The remainder of this paper is structured as follows: An overview of past effort on OLLPs is prvided in Section~\ref{lit_rev}. Section~\ref{problem_defenition} delineates the specific attributes of the open layout problem tackled in this study. In Section~\ref{problem_formulation}, we introduce the mathematical model developed to address this problem. Section~\ref{solution_method} describes the innovative solution approach we devised, detailing the techniques and strategies employed. Section~\ref{exp_setup} outlines the experimental design, including the benchmark problems and methodologies utilized. The results are presented and critically compared with existing literature in Section~\ref{experimental_results}. Section~\ref{discussion} offers additional insights and reflections derived from the study’s findings and lessons learned. Finally, Section~\ref{conclusion} concludes the paper and discusses potential directions for future inquiry.

\section{OLLP: An Overview}\label{lit_rev}
Strategic layout design has shown remarkable potential for cost savings and efficiency gains. For instance, optimal machinery placement in manufacturing can reduce production lead times, material handling costs, and overall operational expenses by up to 35\% \citep{KlaRai2019OptimalFL, alduaij2020adopting}. Similarly, efficient warehouse layouts in logistics and supply chain management improve inventory handling, accelerate order processing, and boost overall performance \citep{diaz2024variable}. Beyond industrial applications, well-optimized layouts in healthcare and retail improve service delivery by enhancing flow and space utilization, as seen in hospitals and shopping malls \citep{meng2021impact}. In office environments, optimized layouts promote increased employee productivity and collaboration \citep{bukTzu2014new, dorrah2021integrated}.

OLLPs represents a significant extension of layout optimization, focusing on arranging rectangular entities on a plane without overlap to minimize transportation costs \citep{perez2021facility}. In OLLPs, entities can be placed anywhere on the plane as long as they meet the non-overlapping constraints. Unlike closed-loop systems with predefined and structured flow paths, open-loop configurations are characterised by unrestricted and flexible pathways. However, this flexibility introduces unique optimisation challenges, making OLLPs more complex to model and solve effectively \citep{niroomand2013mixed, Chae01072006}. The inherent complexity of OLLPs comes from the need to model obstacles, distance metrics, and spatial constraints when planning the paths between entities.
These factors contribute to the classification of OLLP as NP-hard, where the exponential growth in the number of possible configurations renders the problem computationally intractable for exact solutions. Consequently, approximation algorithms and heuristic methods are essential for addressing real-world applications of OLLP \citep{drira2006survey}. 
The evolution of layout optimisation techniques reflects a gradual shift from simplistic models to sophisticated, real-world applications. \citep{kundu2012metaheuristic}.

Early approaches that leveraged the concept of weighted centroids—as introduced in \citep{muther1973systematic} through analyses of material flow and departmental placement—laid the foundation for later centroid-based distance models, offering computationally efficient mechanisms for solving layout problems. However, these models have limited applicability as they overlook practical constraints such as obstacles and non-linear flow paths \citep{MelBoz1996PR}. These limitations were particularly evident in scenarios involving irregularly shaped entities or dynamic flow interactions. The introduction of rectangular layout problems in the late 1970s marked a significant advancement in addressing challenges posed by irregular shapes \citep{thornton1979rectangular}. Since then, a range of mathematical models, along with heuristic and metaheuristic approaches, have been proposed. However, these methods are predominantly based on distance approximations, limiting their accuracy and practical applicability due to inherent suboptimality (see Figure~\ref{fig:megoldas} for an example).

\begin{figure}[ht]
    \centering
       \begin{tikzpicture}[scale=.7]
        \fill[gray!20] (-6.3,4) rectangle (-1.6,1.8);
        \begin{scope}[xshift=-6cm, yshift=3.5cm]
            \draw[thick] (0,0) rectangle (1,0.5); 
            \node[anchor=west] at (1.4,0.25) {Rectangle};
            \draw[line width=2pt,black,dashed](0,-0.5) -- (1,-0.5); 
            \node[anchor=west] at (1.4,-0.5) {Flow Path};
            \draw[thick] (0,-1) -- ++(0.4,0); 
            \draw[thick] (0,-1.3) -- ++(0.4,0); 
            \node[anchor=west] at (1.4,-1.1) {Door};
        \end{scope}
        
        \begin{scope}[xshift=0cm, yshift = -0]
            \draw[thick] (0,0-0.1) rectangle (3,1-0.1);
            \draw[thick] (1.3,1) -- ++(0,-0.2-0.1); 
            \draw[thick] (1.7,1) -- ++(0,-0.2-0.1); 
            
            \draw[thick] (0.5-0.1,1) rectangle (1.5-0.1,4);
            \draw[thick] (1.5-0.1,2.3) -- ++(-0.2,0); 
            \draw[thick] (1.5-0.1,2.7) -- ++(-0.2,0); 
            
            \draw[thick] (1.5+0.1,1) rectangle (2.5+0.1,4);
            \draw[thick] (1.5+0.1,2.3) -- ++(0.2,0); 
            \draw[thick] (1.5+0.1,2.7) -- ++(0.2,0); 
            
            \draw[line width=2pt,black,dashed] (1.5,2.5) -- (1.5,1); 
            
            \node[below] at (1.5,-0.1) {Optimal Solution};
        \end{scope}
        
        \begin{scope}[xshift=6cm]
            \draw[thick] (0-0.1,0) rectangle (1-0.1,3);
            \draw[thick] (1-0.1,1.2) -- ++(-0.2,0); 
            \draw[thick] (1-0.1,1.6) -- ++(-0.2,0); 
            
            \draw[thick] (1,0) rectangle (2,3);
            \draw[thick] (1,1.2) -- ++(0.2,0); 
            \draw[thick] (1,1.6) -- ++(0.2,0); 
            
            \draw[thick] (2+0.1,0) rectangle (3+0.1,3);
            \draw[thick] (2+0.1,1.2) -- ++(0.2,0); 
            \draw[thick] (2+0.1,1.6) -- ++(0.2,0); 
            
            \draw[line width=2pt,black,dashed] (1,1.5) -- (1,3);  
            \draw[line width=2pt,black,dashed](1,3) -- (2,3);  
            \draw[line width=2pt,black,dashed](2,3) -- (2,1.5); 
            
            \node[below] at (1.5,-0.1) {Suboptimal Solution};
        \end{scope}
    \end{tikzpicture}
    \caption{Suboptimality in centroid-based methods}
    \label{fig:megoldas}
\end{figure}

For instance, \citet{niroomand2013mixed} introduced mixed-integer programming models, incorporating rigorous constraints to account for feasible flow paths. However, their reliance on approximated $l_1$ distance metrics limited efficiency in scenarios where optimal paths deviated from linear assumptions. Similarly, \citet{ma2025enhancing} employed a centroid-based approximation to improve urban emergency response, raising concerns about the accuracy of the model. While using $l_1$ distance approximations may be justified in large-scale problems where computing exact distances is computationally prohibitive, this limitation becomes particularly critical in facility layout contexts. Metaheuristic methods are no exception, and their solution quality is influenced by the reliance on approximate distance calculations. For example, \citet{Hauser01112006} and \citet{ONUT2008783} demonstrated the effectiveness of GA and PSO, respectively, in tackling high-dimensional and non-linear optimization challenges. To enhance efficiency and adaptability, many studies have explored the integration of machine learning techniques and hybrid approaches within metaheuristic frameworks.  \citet{powell2010merging} demonstrated the potential of combining artificial intelligence with traditional optimization strategies to address complexity and dynamic constraints. Hybrid approaches, such as simulation-based optimization \citep{YanLiuXu2023simulation} and surrogate models like Kriging \citep{zhang2017novel}, capture complex interactions and enhance computational efficiency while maintaining solution quality. However, the reliance on approximate distance metrics remains a critical limitation, often resulting in suboptimal outcomes. Precise distance calculations become even more complex when the concept of doors or input/output (IO) points is introduced, as it necessitates door-to-door path planning \citep{khalilabadi2024exploiting}.

In layout optimization, door-to-door distance refers to the distance between entry points (doors) of spatial entities such as rooms, departments, or workstations. This metric is critical for accurately modeling flow dynamics within a facility, as it captures the true spatial relationships and movement patterns between different areas \citep{khalilabadi2024exploiting}. By considering the distance between doors, along with traffic volume and directionality, realistic simulations of flow efficiency can be developed, accounting for congestion, bottlenecks, and pathway effectiveness \citep{benson1997doorfast}. In warehouses, for example, reducing door-to-door distances improves efficiency by cutting travel time and reducing congestion along transportation routes \citep{kim2015optimal}.

In the following sections, we present a concise description of the rectangular OLLP with door-to-door distance requirements. For the first time, an exact mathematical model is formulated to address this problem, complemented by a novel solution approach that integrates exact door-to-door path planning with advanced metaheuristic algorithms.

\section{Problem Definition}\label{problem_defenition}
\begin{figure}[ht]
  \centering
  \begin{subfigure}{0.3\textwidth}
    \centering
    \begin{tikzpicture}[scale=0.8, clip]

      \fill[gray!20] (-1.4,-1.95) rectangle (4.4,0.75);
      %
      \draw[thick] (-1,0) rectangle ++(1,0.5);
      \node[anchor=west] at (0.2,0.25) {Rectangle};
      %
      \fill[pattern=north west lines] (-1,-0.6) rectangle ++(1,0.5); 
      \node[anchor=west] at (0.2,-0.34) {Overlapping Area};
      %
      \draw[thick, dashed] (-1,-0.8) -- ++(1,0);
      \node[anchor=west] at (0.2,-0.8) {Flow Path};
      %
      \draw[thick, dashed, red, line width=1.5pt] (-0.8,-1.7) circle [radius=0.2];
      \node[anchor=west] at (0.2,-1.7) {Intersection};
      %
      \draw[thick] (-1,-1.3) -- ++(0.4,0); 
      \draw[thick] (-1,-1.1) -- ++(0.4,0); 
      \node[anchor=west] at (0.2,-1.2) {Door};
    \end{tikzpicture}
    \label{fig:legend_section}
  \end{subfigure}
  \begin{subfigure}{0.3\textwidth}
    \centering
    \begin{tikzpicture}[scale=0.8]
      \def\rectwidth{1}
      \def\rectheight{3}
      \draw[thick] (0,0) rectangle ++(\rectwidth,\rectheight)
            node[midway, rotate=90, anchor=center, yshift=0.75cm] {Rectangle 1};
      \draw[thick] ($(0,0) + (0,\rectheight/2 - 0.2)$) -- ++(0.2,0);
      \draw[thick] ($(0,0) + (0,\rectheight/2 + 0.2)$) -- ++(0.2,0);
      \draw[thick] (0.5,0.5) rectangle ++(\rectwidth,\rectheight)
            node[midway, rotate=90, anchor=center, yshift=-0.85cm] {Rectangle 2};
      \draw[thick] ($(0.5,0.5) + (\rectwidth, \rectheight/2 - 0.2)$) -- ++(-0.2,0);
      \draw[thick] ($(0.5,0.5) + (\rectwidth, \rectheight/2 + 0.2)$) -- ++(-0.2,0);
      \fill[pattern=north west lines] (0.5,0.5) rectangle ++(0.5,2.5);
    \end{tikzpicture}
    \caption{Overlapping}
    \label{fig:overlap_boxes}
  \end{subfigure}
  \begin{subfigure}{0.3\textwidth}
    \centering
    \begin{tikzpicture}[scale=0.8]
      \begin{scope}[shift={(-2cm,0)}]
        \def\rectwidth{1}
        \def\rectheight{3}
        \draw[thick] (0,0) rectangle ++(\rectwidth,\rectheight)
             node[midway, rotate=90, anchor=center] {Rectangle 1};
        \draw[thick] ($(0,0) + (\rectwidth,\rectheight/2 - 0.2)$) -- ++(-0.2,0);
        \draw[thick] ($(0,0) + (\rectwidth,\rectheight/2 + 0.2)$) -- ++(-0.2,0);
        \draw[thick] (4,0) rectangle ++(\rectwidth,\rectheight)
             node[midway, rotate=90, anchor=center] {Rectangle 2};
        \draw[thick] ($(4,0) + (0,\rectheight/2 - 0.2)$) -- ++(0.2,0);
        \draw[thick] ($(4,0) + (0,\rectheight/2 + 0.2)$) -- ++(0.2,0);
        \draw[thick, dashed] (0,\rectheight) -- (5,0);
        \draw[thick, dashed, red, line width=1.5pt] (\rectwidth, \rectheight-\rectwidth+0.36) circle [radius=0.15];
        \draw[thick, dashed, red, line width=1.5pt] (4, \rectwidth-0.38) circle [radius=0.15];
      \end{scope}
    \end{tikzpicture}
    \caption{Intersection}
    \label{fig:intersection_boxes}
  \end{subfigure}
  
  \caption{Illustrations of (a)Overlapping and, (b)Intersection of a path and the edges.}
  \label{fig:layout_figures}
\end{figure}
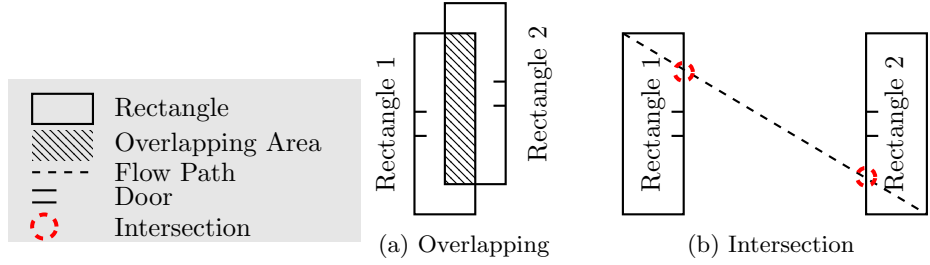

We consider a specific type of layout problem that involves placing a number of rectangles of different dimensions on a plane without overlapping (Figure~\ref{fig:overlap_boxes}). These rectangles represent different objects across different domains. A typical example is the arrangement of machines or cells in a factory, where the objective is to minimise the total transportation cost of semi-finished products \citep{shi2022study}.

To provide a more realistic representation of the problem, we assume each rectangle (cell) is equipped with a pick-up point positioned at the midpoint of one of its edges. Transportation occurs only between these pick-up points, here called "doors". A general constraint is that the rectangles must have vertical and horizontal edges and can be rotated by 90, 180, or 270 degrees (\(\rightarrow\) Section \ref{problem_formulation} \(\rightarrow\) Figure~\ref{fig:rotation_box}), only. Transportation paths cannot cross the rectangles, and any intersection (Figure \ref{fig:intersection_boxes}) is infeasible; they must move either outside or along the edges of the rectangles.

\section{Problem Formulation}\label{problem_formulation}
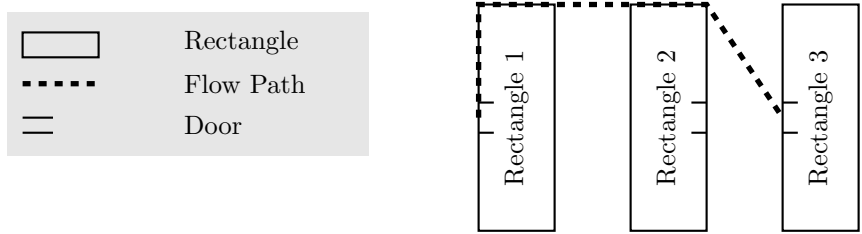
\begin{figure}[h]
    \centering

    \begin{subfigure}[b]{0.8\textwidth}

\begin{tikzpicture}[scale=1]
    \def\rectwidth{1}
    \def\rectheight{3}
    
    \coordinate (Rect1) at (2,0);
    \coordinate (Rect2) at (4,0);
    \coordinate (Rect3) at (6,0);
    
    \draw[thick] (Rect1) rectangle ++(\rectwidth,\rectheight) 
        node[midway, rotate=90, anchor=center] {Rectangle 1};
    \draw[thick] (Rect2) rectangle ++(\rectwidth,\rectheight) 
        node[midway, rotate=90, anchor=center] {Rectangle 2};
    \draw[thick] (Rect3) rectangle ++(\rectwidth,\rectheight) 
        node[midway, rotate=90, anchor=center] {Rectangle 3};
    
    \draw[thick] ($(Rect1) + (0,\rectheight/2 - 0.2)$) -- ++(0.2,0);
    \draw[thick] ($(Rect1) + (0,\rectheight/2 + 0.2)$) -- ++(0.2,0);

    \draw[thick] ($(Rect2) + (\rectwidth,\rectheight/2 - 0.2)$) -- ++(-0.2,0);
    \draw[thick] ($(Rect2) + (\rectwidth,\rectheight/2 + 0.2)$) -- ++(-0.2,0);

    \draw[thick] ($(Rect3) + (0,\rectheight/2 - 0.2)$) -- ++(0.2,0);
    \draw[thick] ($(Rect3) + (0,\rectheight/2 + 0.2)$) -- ++(0.2,0);
    
    \draw[line width=2pt,black,dashed] 
        ($(Rect1) + (0,\rectheight/2)$) -- ++(0,\rectheight/2) -- ++(\rectwidth,0) -- ($(Rect2) + (\rectwidth, \rectheight)$) -- ($(Rect3) + (0,\rectheight/2)$);
    
    \begin{scope}[shift={(-4,2.3)}] 
        \fill[gray!20] (-0.2,-1.3) rectangle (4.55,.6);
        
        \draw[thick] (0,0) rectangle ++(1.,0.33);
        \node[anchor=west] at (2,0.2) {Rectangle};

        \draw[line width=2pt,black,dashed] (0,-0.35) -- ++(1.,0);
        \node[anchor=west] at (2,-0.35) {Flow Path};

        \draw[thick] (0,-1.) -- ++(0.4,0);
        \draw[thick] (0,-.8) -- ++(0.4,0);
        \node[anchor=west] at (2,-.9) {Door};
    \end{scope}
\end{tikzpicture}

    \end{subfigure}
    
    \caption{Exact feasible $l_2$ distances.}
    \label{fig:l1_l2_distance_layout}
\end{figure}

An exact mathematical model is provided for the case where distances are measured using the Euclidean metric. It should be noted that  in mathematics, there are infinitely many other ways to measure distances, commonly represented as $l_p$ norms, where $1 \leq p \leq \infty$. The Euclidean distance corresponds to the $l_2$ norm. An example of a  feasible \( l_2 \) path is illustrated in Figure~\ref{fig:l1_l2_distance_layout}.
In engineering, two other distance measures are particularly significant: the $l_1$ and $l_{\infty}$ norms. The $l_1$ distance, also known as the Manhattan distance, is the sum of the absolute differences between the coordinates. The $l_{\infty}$ distance, on the other hand, is defined as the maximum of the absolute differences between coordinates. Each of these distance metrics requires a distinct mathematical model to determine the path that is optimal (i.e., minimal) with respect to the given distance measure and thus are not discussed here.

\subsection{Mathematical Model for OLLPs with Exact Euclidean Distances}
\begin{table}[ht]
    \centering
    \caption{Sets, Parameters and Variables }
    \label{Nomenclature}
    \begin{tabularx}{\textwidth}{@{}lX@{}}
        \toprule
        \textbf{Symbol} & \textbf{Description} \\
        \midrule
        $I$ & Set of cell (rectangle) indices.  \textbf{Alias:} \(i,\, j,\, p\) \\
        $K=\{1,2,3,4,5\}$ & Set of node (vertices + door) indices for each cell\\
        \(R=\{1,2,3,4\}\) & Edge index\\
        \(Q=\{1,2,3,4\}\) & Orientation index\\
        \midrule
         \(s_i (t_i)\) & Length of the vertical (horizontal) edge of cell \(i\) \\
          \(F_{i,j}\) & The flow between rectangles $i$ and $j$ \\
          \(\mathbf{A}\) & Node–arc incidence matrix of the network connecting door (pick–up) points. \\
        \midrule
        \(a_i\) (\(b_i\)) & Horizontal (vertical) coordinate of the center of cell \(i\) \\
        \(e_{ij}, f_{ij}\) (\(g_{ij}, h_{ij}\)) & Horizontal (vertical) distances between the centers of cells \(i\) and \(j\) \\
        \(z_i\) & Binary variable; 1 if cell \(i\) is in a vertical position, 0 if horizontal \\
        \(x_{i_k}\) (\(y_{i_k}\)) & Coordinate of node \(k\in K\) of cell \(i\)\\
        \(\lambda_{i_q}\) & Binary variable describing the door position for cell \(i\) \\
       
        \(d_{i_kj_l}\) & Distance between nodes \(i_k\) and \(j_l\) \\
        \(M_{i_kj_lp_r}\) & Penalty if the line connecting \(i_k\) to \(j_l\) crosses the edge \(r\) of cell \(p\) \\
        \(\gamma_{i_kj_lp_r},\, \delta_{i_kj_lp_r}\) & Auxiliary variables for penalty calculation \\
        \(u_{i_kj_lp_r},\, v_{i_kj_lp_r},\, w_{i_kj_lp_r},\, z_{i_kj_lp_r},\, m_{i_kj_lp_r}\) & Binary auxiliary variables for penalty calculation \\
        \(\alpha_{ij}, \beta_{ij}, \mu_{ij}\) & Binary auxiliary variables for non-overlapping conditions between cells \(i\) and \(j\) \\
        \(\mathbf{e}_{i_5}\) & Unit vector associated with cell \(i\) (at its pick-up point) \\
        \(\mathbf{g}_{ij}\) & Path vector between the pick-up points \(i_5\) and \(j_5\) \\
        \bottomrule
    \end{tabularx}
\end{table}

The exact description of the objective function requires certain data from the constraint set to ensure the non-overlapping property of the layout. \citet{mont} was the first in discussing this type of constraint. More complicated models can be found in \citep{hk, tompkins1996facilities, meller, sherali}, and \citep{das}. Any of the aforementioned sets of non-overlapping constraints can be used in an exact model, although significant differences may occur in the effectiveness of the model. Here we use the latter.
In what follows, we will outline these basic constraints along with specific constraints required for open layout problems for $l_2$ distances considering the positions of the doors which are unique to this paper. The nomenclature is defined in Table~\ref{Nomenclature}.

\begin{figure}[ht]
    \centering
    \begin{tikzpicture}[scale=0.8]
        \def\boxwidth{1}
        \def\boxheight{2.5} 
        \def\rowgap{1.5} 
        \begin{scope}[shift={(9.5,0)}]
            \draw[thick] (0,0) rectangle ++(\boxwidth,\boxheight);
            \draw[thick] (0,0+\boxheight/2 - 0.2) -- ++(0.2,0); 
            \draw[thick] (0,0+\boxheight/2 + 0.2) -- ++(0.2,0); 
            \node[anchor=north] at (0.5, -0.2) {$\lambda_{i_4} = 1$};
            \node[anchor=north] at (0.5, -0.9) {$\sum\limits_{j=1,2,3}(\lambda_{i_j}) = 0$};

            \node[anchor=north] at (0.5, \boxheight+0.7) {$s_i$}; 
            \node[anchor=east] at (-0.2, \boxheight/2) {$t_i$};  
            \node[anchor=north] at (-0.2, \boxheight+0.7) {$i_1$}; 
            \node[anchor=north] at (1.2, \boxheight+0.7) {$i_2$}; 
            \node[anchor=north] at (-0.2,0.35) {$i_4$}; 
            \node[anchor=north] at (1.2,0.35) {$i_3$}; %
            
        \end{scope}

        \begin{scope}[shift={(1.4,0)}]
            \draw[thick] (3,0) rectangle ++(\boxheight,\boxwidth); 
            \draw[thick] (3+\boxheight/2 - 0.2,0+\boxwidth) -- ++(0,-0.2); 
            \draw[thick] (3+\boxheight/2 + 0.2,0+\boxwidth) -- ++(0,-0.2); 

            \node[anchor=north] at (4.5, -0.2) {$\lambda_{i_3} = 1$};
            \node[anchor=north] at (4.5, -.9) {$\sum\limits_{j=1,2,4}(\lambda_{i_j}) = 0$};
            \node[anchor=north] at (2.5, \boxwidth-.25) {$s_i$}; 
            \node[anchor=east] at (4.5, \boxheight/2) {$t_i$};  
            \node[anchor=north] at (3-0.2, \boxwidth+0.7) {$i_1$}; 
            \node[anchor=north] at (5.5+.2, \boxwidth+0.7) {$i_2$}; 
            \node[anchor=north] at (3-0.2,0.35) {$i_4$}; 
            \node[anchor=north] at (5.5+.2,0.35) {$i_3$};
        \end{scope}
        \begin{scope}[shift={(-4.9,0)}]
            \draw[thick] (6,0) rectangle ++(\boxwidth,\boxheight);
            \draw[thick] (6+\boxwidth,0+\boxheight/2 - 0.2) -- ++(-0.2,0); 
            \draw[thick] (6+\boxwidth,0+\boxheight/2 + 0.2) -- ++(-0.2,0); 
            \node[anchor=north] at (6.5, -.2) {$\lambda_{i_2} = 1$};
            \node[anchor=north] at (6.5, -.9) {$\sum\limits_{j=1,3,4}(\lambda_{i_j}) = 0$};
            \node[anchor=north] at (6.5, \boxheight+.75) {$s_i$}; 
            \node[anchor=east] at (6, \boxheight/2) {$t_i$};  
        \end{scope}

        \begin{scope}[shift={(-14.4,0)}]
            \draw[thick] (10.5,0) rectangle ++(\boxheight,\boxwidth);
            \draw[thick] (10.5+\boxheight/2 - 0.2,0) -- ++(0,0.2); 
            \draw[thick] (10.5+\boxheight/2 + 0.2,0) -- ++(0,0.2); 

            \node[anchor=north] at (12, -.2) {$\lambda_{i_1} = 1$};
            \node[anchor=north] at (12, -.9) {$\sum\limits_{j=2,3,4}(\lambda_{i_j}) = 0$};
            \node[anchor=north] at (11.75, \boxwidth+.75) {$t_i$}; 
            \node[anchor=east] at (10.5, \boxheight/2-.75) {$s_i$};  
            \node[anchor=north] at (10.5-0.2, \boxwidth+0.7) {$i_1$}; 
            \node[anchor=north] at (13.5-.2, \boxwidth+0.7) {$i_2$}; 
            \node[anchor=north] at (10.5-0.2,0.35) {$i_4$}; 
            \node[anchor=north] at (13.5-.2,0.35) {$i_3$};
        \end{scope}


    \end{tikzpicture}
    \caption{Orientation of the $i^{th}$ rectangular box ($s_i$,$t_i$), described by \(\lambda_{i_j}\).}
    \label{fig:rotation_box}
\end{figure}
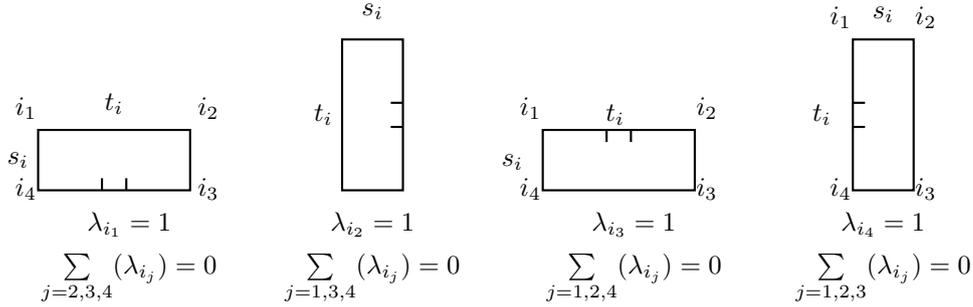
Each cell contains five key points: its four corner points and a pick-up point. When the cell is rotated, the locations and indices of these points are illustrated in Figure \ref{fig:rotation_box}. The upper-left corner is always assigned index $1$, while the upper-right corner is consistently indexed as $2$, and so on. However, the position of the pick-up point varies depending on the cell's rotation. In its default position ($0 ^{\circ}$  rotation), the length of the vertical (horizontal) edges of cell $i$ is $s_i$($t_i$). Without loss of generality, we assume that the door (or the edge containing the door) is located below the centroid in the default position. The coordinates of the doors of the \(i^{th}\) rectangle depend only on its orientation/rotation on the plane. The orientation of  the \(i^{th}\) cell  is described by four binary variables, \(\lambda_{i_j};j\in\{1,2,3,4\}\). When \(\lambda_{i_1}=1\) the door is located below the centroid of the rectangle $i$, when \(\lambda_{i_2}=1\) the door is located to the right of the centroid, when \(\lambda_{i_3}=1\) the door is located above the centroid, and finally \(\lambda_{i_4}=1\)  means the door is located to the left of the centroid. \(\lambda_{i_j}=1\) when the pick-up point is in the specified position dictated by $j$, otherwise $0$ such that $\Sigma_{j=1}^4{\lambda_{ij}}=1 \textbf{  }\forall i$. The complete mathematical model will be described as follows:
\allowdisplaybreaks
{\small
\begin{flalign}
   & \min \sum_{i,j} F_{ij}\, {\bf d^T}_{ij} {\bf g}_{ij} && && \label{obj} \\
   & {\bf A}\, {\bf g}_{ij} = {\bf e}_{i_5} - {\bf e}_{j_5} && \forall\, i,j\in I && \label{path} \\
   & e_{ij}+f_{ij}-\frac{1-z_i}{2}t_i-\frac{z_i}{2}s_i-\frac{1-z_j}{2}t_j-\frac{z_j}{2}s_j \geq -M\mu_{ij} && \forall\, i,j\in I && \label{nono1}\\
   & g_{ij}\!+\!h_{ij}-\!\frac{1\!-\!z_i}{2}s_i\!-\!\frac{z_i}{2}t_i\!-\!\frac{1\!-\!z_j}{2}s_j\!-\!\frac{z_j}{2}t_j\! \geq\! -M(1\!-\!\mu_{ij}) && \forall\, i,j\in I && \label{nono2}\\
   & a_i-a_j = e_{ij}-f_{ij} && \forall\, i,j\in I && \label{nono3}\\
   & b_i-b_j = g_{ij}-h_{ij} && \forall\, i,j\in I && \label{nono4}\\
   & e_{ij}\leq M\alpha_{ij} && \forall\, i,j\in I && \label{nono5}\\
   & f_{ij}\leq M(1-\alpha_{ij}) && \forall\, i,j\in I && \label{nono6}\\
   & g_{ij}\leq M\beta_{ij} && \forall\, i,j\in I && \label{nono7}\\
   & h_{ij}\leq M(1-\beta_{ij}) && \forall\, i,j\in I && \label{nono8}\\[1mm]
   & x_{i_1} = x_{i_4} = a_i - \frac{1-z_i}{2}t_i-\frac{z_i}{2}s_i && \forall\, i\in I && \label{e31} \\
   & x_{i_2} = x_{i_3} = a_i + \frac{1-z_i}{2}t_i+\frac{z_i}{2}s_i && \forall\, i\in I && \label{e32} \\
   & y_{i_1} = y_{i_2} = b_i + \frac{1-z_i}{2}s_i+\frac{z_i}{2}t_i && \forall\, i\in I && \label{e33} \\
   & y_{i_3} = y_{i_4} = b_i - \frac{1-z_i}{2}s_i-\frac{z_i}{2}t_i && \forall\, i\in I && \label{e34} \\
   & z_i = \lambda_{i_2}+\lambda_{i_4} && \forall\, i\in I && \label{e09} \\
   & 1-z_i = \lambda_{i_1}+\lambda_{i_3} && \forall\, i\in I && \label{e10} \\
   & x_{i_5} = a_i+(\lambda_{i_2}-\lambda_{i_4})\frac{s_i}{2} && \forall\, i\in I && \label{e11} \\
   & y_{i_5} = b_i+(\lambda_{i_3}-\lambda_{i_1})\frac{s_i}{2} && \forall\, i\in I && \label{e12} \\[1mm]
   & d_{i_1i_2} = d_{i_3i_4} = z_i s_i+(1-z_i)t_i && \forall\, i\in I && \label{adj_dis1} \\
   & d_{i_2i_3} = d_{i_1i_4} = z_i t_i+(1-z_i)s_i && \forall\, i\in I && \label{adj_dis2} \\
   & d_{i_1i_3} = d_{i_2i_4} = M && \forall\, i\in I && \label{opposite_cells} \\[1mm]
   & d_{i_5i_1} = (\lambda_{i_1}+\lambda_{i_2})M+(\lambda_{i_3}+\lambda_{i_4})\frac{t_i}{2} && \forall\, i\in I && \label{pick1} \\
   & d_{i_5i_2} = (\lambda_{i_1}+\lambda_{i_4})M+(\lambda_{i_2}+\lambda_{i_3})\frac{t_i}{2} && \forall\, i\in I && \label{pick2} \\
   & d_{i_5i_3} = (\lambda_{i_3}+\lambda_{i_4})M+(\lambda_{i_1}+\lambda_{i_2})\frac{t_i}{2} && \forall\, i\in I && \label{pick3} \\
   & d_{i_5i_4} = (\lambda_{i_2}+\lambda_{i_3})M+(\lambda_{i_1}+\lambda_{i_4})\frac{t_i}{2} && \forall\, i\in I && \label{pick4} \\[1mm]
   & d_{i_kj_l} = \sqrt{(x_{i_k}-x_{j_l})^2+(y_{i_k}-y_{j_l})^2} + \sum_{p=1}^n \sum_{r=1}^4 M_{i_kj_lp_r} && \forall\, i,j\in I;\; k,l\in K && \label{2cell20} \\[1mm]
   & x_{p_2} = x_{p_3} = \gamma_{i_kj_lp_1}x_{i_k} + (1-\gamma_{i_kj_lp_1})x_{j_l} && \forall\, i,j,p\in I;\; k,l\in K && \label{2cell21} \\
   & \gamma_{i_kj_lp_1}y_{i_k}\! + \!(1\!-\!\gamma_{i_kj_lp_1})y_{j_l} \!= \!\delta_{i_kj_lp_1}y_{p_2} \!+\! (1\!-\!\delta_{i_kj_lp_1})y_{p_3} && \forall\, i,j,p\in I;\; k,l\in K && \label{2cell22} \\
   & x_{p_1} = x_{p_4} = \gamma_{i_kj_lp_3}x_{i_k} + (1-\gamma_{i_kj_lp_3})x_{j_l} && \forall\, i,j,p\in I;\; k,l\in K && \label{2cell23} \\
   & \gamma_{i_kj_lp_3}y_{i_k} \!+\! (1\!-\!\gamma_{i_kj_lp_3})y_{j_l} \!= \!\delta_{i_kj_lp_3}y_{p_1}\! +\! (1\!-\!\delta_{i_kj_lp_3})y_{p_4} && \forall\, i,j,p\in I;\; k,l\in K && \label{2cell24} \\
   & y_{p_4} = y_{p_3} = \delta_{i_kj_lp_2}y_{i_k} + (1-\delta_{i_kj_lp_2})y_{j_l} && \forall\, i,j,p\in I;\; k,l\in K && \label{2cell25} \\
   & \delta_{i_kj_lp_2}x_{i_k} \!+\! (1\!-\!\delta_{i_kj_lp_2})x_{j_l} \!= \!\gamma_{i_kj_lp_2}x_{p_4} \!+ \!(1\!-\!\gamma_{i_kj_lp_2})x_{p_3} && \forall\, i,j,p\in I;\; k,l\in K&& \label{2cell26} \\
   & y_{p_1} = y_{p_2} = \delta_{i_kj_lp_4}y_{i_k} + (1-\delta_{i_kj_lp_4})y_{j_l} && \forall\, i,j,p\in I;\; k,l\in K && \label{2cell27} \\
   & \delta_{i_kj_lp_4}x_{i_k} \!+ \!(1\!-\!\delta_{i_kj_lp_4})x_{j_l} \!= \!\gamma_{i_kj_lp_4}x_{p_1} \!+ \!(1\!-\!\gamma_{i_kj_lp_4})x_{p_2} && \forall\, i,j,p\in I;\; k,l\in K && \label{2cell28} \\[1mm]
   & 1-\gamma_{i_kj_lp_r}-u_{i_kj_lp_r}M\leq 0 && \forall\, \!i,j,p\!\in \!I;\;\! k,l\!\in\! K;\! r\!\in\! R&& \label{2cell29} \\
   & \gamma_{i_kj_lp_r}-v_{i_kj_lp_r}M\leq 0 &&\forall\,\! i,j,p\!\in \!I;\;\! k,l\!\in \!K; \!r\!\in \!R && \label{2cell30} \\
   & 1-\delta_{i_kj_lp_r}-w_{i_kj_lp_r}M\leq 0 && \forall\,\! i,j,p\!\in \!I;\;\! k,l\!\in\! K; \!r\!\in\! R && \label{2cell31} \\
   & \delta_{i_kj_lp_r}-z_{i_kj_lp_r}M\leq 0 && \forall\, \!i,j,p\!\in\! I;\;\! k,l\!\in\! K; r\!\in\! R && \label{2cell32} \\
   & u_{i_kj_lp_r}+v_{i_kj_lp_r}+w_{i_kj_lp_r}+z_{i_kj_lp_r}\leq 3+m_{i_kj_lp_r} &&\forall\, \!i,j,p\!\in\! I;\; \!k,l\!\in\! K;\! r\!\in\! R && \label{2cell33} \\
   & M_{i_kj_lp_r} = m_{i_kj_lp_r}\, M && \forall\, \!i,j,p\!\in \!I;\;\! k,l\!\in\! K;\! r\!\in\! R && \label{2cell34} \\[1mm]
   & e_{ij},\, f_{ij},\, g_{ij},\, h_{ij} \geq 0 && \forall\, i,j\in I && \label{e170}\\
   & d_{i_kj_l}, g_{i_kj_l} \geq 0 && \forall\, i,j\in I,\; \forall\, k,l\in K && \label{e171}\\
   & u_{i_kj_lp_r},\, v_{i_kj_lp_r},\, w_{i_kj_lp_r},\, z_{i_kj_lp_r},\, m_{i_kj_lp_r}\in \{0,1\} && \forall\, \!i,j,p\!\in \!I;\;\! k,l\!\in\! K; \!r\!\in\! R&& \label{e173}\\
   & \mu_{ij},\, \alpha_{ij},\, \beta_{ij},\, z_i,\, \lambda_{i_q} \in \{0,1\}&& \forall\, i,j\in I;\; q \in Q&& \label{e174}
\end{flalign}}

The objective function (\ref{obj}) seeks to minimize the total flow between cells, weighted by the distances between their respective pick-up points. These pick-up points represent the vertices of a graph, with edges formed between vertices whenever a direct movement is feasible. Consequently, the edge set dynamically depends on the spatial arrangement of the rectangles. The lengths of these edges correspond to the $l_2$ distances between vertices, which are inherently influenced by the locations of the rectangles.

The constraint \ref{path} establishes a path from pick-up point $i_5$ to pick-up point $j_5$ for all non-negative flow values $f_{ij}$, where ${\bf A}$ is the node-arc adjacency matrix of the complete graph, consisting of the critical points of all cells. Constraints \ref{nono1} - \ref{nono8} guarantee the non-overlapping of cells $i$ and $j$. Two cells are overlapping if and only if their centers are too close to each other. The minimal required horizontal (vertical) distance such that two cells are not
overlapping is half the sum of the length of their edges in the horizontal (vertical) position. The sum depends on the rotation of the cells. Let $e_{ij}=\max\{0,x_i-x_j\}$, $f_{ij}=\max\{0,x_j-x_i\}$, $g_{ij}=\max\{0,y_i-y_j\}$ and $h_{ij}=\max\{0,y_j-y_i\}$. Notice that $e_{ij}+f_{ij}$ $(g_{ij}+h_{ij})$ is the horizontal(vertical) distance of the centers of the cells $i$ and $j$. If there is no horizontal (vertical) overlap, then the distance must be at least as long as the sum of the two horizontal (vertical) half edges. This requirement is described by \ref{nono1} and \ref{nono2}. If the binary variable $\mu_{ij}$ is 1, then cells $i$ and $j$ are not overlapping vertically,  and if it is 0, then cells $i$ and $j$ are not overlapping horizontally. It is difficult to use the formulae of $e_{ij}$, $f_{ij}$, $g_{ij}$, and $h_{ij}$ explicitly in an optimization problem; therefore, they are described implicitly by \ref{nono3}-\ref{nono8}. Here $\alpha_{ij}$ is a binary variable: it is 1 if $x_i\geq x_j$; similarly, $\beta_{ij}$ is 1 if $y_i\geq y_j$. Constraints \ref{e31} and \ref{e32} describe the horizontal coordinates of the vertices of cell $i$. Constraints \ref{e33}-\ref{e34} indicate the vertical coordinates of the vertices of cell $i$. Constraints \ref{e09} and \ref{e10} determine if the cell is in horizontal or vertical position. It should be remembered that a cell is in a horizontal position by definition if the pick-up point is below or above the center, and it is in a vertical position if the pick-up point is beside the center. Thus, the $\lambda$ and $z$ variables are not independent. 
Constraints \ref{e11} and \ref{e12} identify the coordinates of the pick-up point of cell $i$. 
In our  model we use the Euclidean distance for the $i_k$ and $j_l$ points if there are no cells between the points. Otherwise the distance will be $M$. Based on equations \ref{adj_dis1}-\ref{adj_dis2} the distance between two adjacent vertices of cell $i$ is equal to the length of the edge. Constraint \ref{opposite_cells} prohibits direct movement between opposite vertices. The constraints \ref{pick1} - \ref{pick4} determines the distance between the pick-up point of cell $i$ and its adjacent vertices: it is half of the length of the edge, i.e., $\frac{t_i}{2}$. Direct movement between the pick-up point and the opposite vertices of the same cell is not feasible; hence, the distance in this case is $M$.

Constraints \ref{2cell20} - \ref{2cell34} determine the distance between points $i_k$ and $j_l$: it is the Euclidean distance if there are no cells between them. Otherwise, the distance is increased by $M$. Here $M_{i_kj_lp_r}$ is a penalty variable; it is $M$ if the line $i_kj_l$ crosses the $r$-th edge of cell $p$
and it is 0 otherwise. If $r=1$, then the $r$-th edge of cell $p$ is the vertical edge between points $p_2$ and $p_3$.
Then there exists the unique value $\gamma_{i_kj_lp_1}$ for which \ref{2cell21} holds. Moreover there exists the value $\delta_{i_kj_lp1}$ which satisfies \ref{2cell22}. If the segment between $i_k$ and $j_l$ crosses the $1^{st}$ edge of cell $p$, then both of the values $\gamma_{i_kj_lp_1}$ and $\delta_{i_kj_lp_1}$ are
between $0$ and $1$, and \ref{2cell21} determines the horizontal, and \ref{2cell22} determines the vertical coordinate of the cross-point. Similarly, if $r=3$, then there exist the same $\gamma_{i_kj_lp_3}$ and $\delta_{i_kj_lp_3}$ for \ref{2cell23} and \ref{2cell24}.
If $r=2$ or $r=4$, then the edge of cell $p$ is horizontal, thus the vertical equation is described first. If $r=2$, then there exist the only $\delta_{i_kj_lp_2}$ and $\gamma_{i_kj_lp_2}$ which satisfy \ref{2cell25} and \ref{2cell26}. Similarly, if $r=4$, then there exist the same $\delta_{i_kj_lp_4}$ and $\gamma_{i_kj_lp_4}$ for \ref{2cell27} and \ref{2cell28}. (For example if $r=4$, then the segment between $i_k$ and $j_l$ crosses the 4th edge of cell $p$ if and only if, then both of the values $\delta_{i_kj_lp_4}$ and $\gamma_{i_kj_lp_4}$ are between 0 and 1.)
If $0<\gamma_{i_kj_lp_r}<1$, then for feasibility of \ref{2cell29} and \ref{2cell30} the values of the two binary variables need to be equal to 1.  Similarly, in \ref{2cell31} and \ref{2cell32} $w$ and $z$ are both equal to 1 if $0<\delta_{i_kj_lp_r}<1$.
If the segment between $i_k$ and $j_l$ crosses the $r$-th edge of cell $p$, then all four binary variables are 1. Therefore in \ref{2cell33} the value of $m$ is 1 if the segment crosses the edge, which leads to the realization of the penalty in \ref{2cell34}.

\section{Methods}\label{solution_method}
Exact methods struggle to scale as the problem size increases due to the inherent NP-hardness and combinatorial complexity of OLLPs \citep{Bennell01102002}. In the case of OLLP with door-to-door distances, even the smallest problem instances require significant time to solve due to the nonlinearity introduced by $l_2$ distance calculations and the large number of binary variables. Consequently, metaheuristic algorithms have emerged as crucial tools for efficiently exploring the solution space and identifying high-quality solutions. Metaheuristics, such as evolutionary algorithms and swarm intelligence, offer a balance between exploration and exploitation, allowing efficient navigation of the vast solution space \citep{TONGUR2020951}.

Metaheuristics rely on solution representation and efficient search operators to explore and refine candidate layouts iteratively. The effectiveness of these algorithms depends on how solutions are encoded and decoded, as these representations directly impact the search process. A well-structured encoding ensures that the evolutionary operators—such as crossover, mutation, and selection—can generate diverse and feasible solutions while preserving spatial constraints. Decoding, in turn, translates these encoded representations into real-world layouts, enabling the evaluation of objective functions and guiding the optimization process. Thus, an effective encoding-decoding mechanism is crucial for leveraging the full potential of metaheuristic algorithms in solving OLLP.
\subsection{Encoding}
\begin{figure}[htbp]
    \centering
    \includegraphics[width=0.9\textwidth]{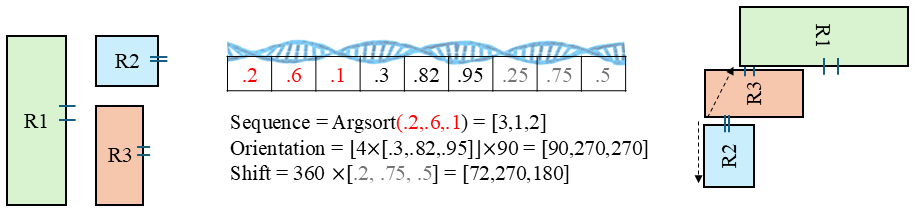}
    \caption{Illustration of a chromosome for a problem with $3$ rectangles.}
    \label{fig:chromosome}
\end{figure}

A feasible layout must be incrementally constructed, ensuring non-overlapping placements, valid door-to-door paths, and optimized spatial arrangements. This process involves determining the insertion sequence of rectangles, their rotation, and their positioning while maintaining feasibility constraints. The proposed encoding comprises three sections, each with a length of $n$, where $n$ denotes the total number of rectangles. The required information for constructing a feasible solution is extracted from the chromosome as described below. A visual example is provided in Figure~\ref{fig:chromosome}.
\begin{enumerate}
    \item \textbf{Insertion Order:} A floating-point representation defining the sequence of rectangle placements. The sequence is sorted based on numerical values, ensuring that lower values correspond to earlier placement. Each rectangle is initially located at the center of the coordinates.
    \item \textbf{Rotation Angle:} The values in this section of the chromosome are mapped to one of $\{0^\circ, 90^\circ, 180^\circ, 270^\circ\}$ by a function $
f(x) = \left\lfloor 4x \right\rfloor \cdot 90^\circ, \quad \text{where } x \in [0, 1]
$, dictating each rectangle’s orientation.
    \item \textbf{Shift Angle:} Each rectangle is associated with a continuous value \(\theta \in [0, 1]\), which dictates the initial displacement direction. The angle along which the rectangle will be shifted is determined by multiplying the value of \(\theta\) by 360. The shifting process continues iteratively in the given direction until no further overlapping with other rectangles exists.
\end{enumerate}

\subsection{Path Computation and Adjacancy Graph Construction}

Once all rectangles are placed without overlaps, an adjacency graph is created to represent accessible node pairs. Each rectangle contributes five key nodes: a door at the center of one edge and one at each vertex.

\textbf{Intra-rectangle connections:}
Feasible paths connect vertices on the same edge, with the path length equal to the edge's length. The door, located at the midpoint of an edge, is connected to its adjacent vertices with paths of half the edge's length.

\textbf{Inter-rectangle connections:}
Feasible paths between nodes on different rectangles are considered if they can be established without intersecting the edges of other rectangles. To decide if the path intersects any edges of positioned rectangles, we use an orientation-based approach \citep{preparata1985computational}.
\begin{definition}[The Orientation Function]
For any three points \(P = (P_x, P_y)\), \(Q = (Q_x, Q_y)\), and \(R = (R_x, R_y)\), we define the \emph{orientation} as:
\[
\operatorname{orient}(P, Q, R) = (Q_x - P_x)(R_y - P_y) - (Q_y - P_y)(R_x - P_x).
\]

\end{definition}

For segments \(AB\) and \(CD\) to intersect in the general case (i.e., when the endpoints are not all collinear), the following conditions must be satisfied:
\[
\operatorname{orient}(A, B, C) \cdot \operatorname{orient}(A, B, D) \leq 0,
\]
and
\[
\operatorname{orient}(C, D, A) \cdot \operatorname{orient}(C, D, B) \leq 0.
\]
When one or more of the orientation values is zero, the corresponding points are collinear. In such instances, it is necessary to verify whether the collinear point lies on the line segment formed by the other two endpoints. For example, if
\[
\operatorname{orient}(A, B, C) = 0,
\]
then point \(C\) is on segment \(AB\) if and only if
\[
\min(x_1, x_2) \leq x_3 \leq \max(x_1, x_2) \quad \text{and} \quad \min(y_1, y_2) \leq y_3 \leq \max(y_1, y_2).
\]

If the path between two points is feasible,
its length is computed using the Euclidean metric. Additionally, paths traversing along shared edges of neighboring rectangles are considered feasible. Two rectangles may share an edge or even a single point (at their edges or vertices) as long as their boundaries do not cross. To enforce this and to mitigate spurious intersections caused by floating-point precision issues, a small buffer threshold is incorporated during the feasibility checks. This ensures computational robustness and the accuracy of path evaluations. This structured approach guarantees that all potential paths are evaluated effectively, minimizing the overall transportation cost while adhering to spatial constraints. \\
With the adjacency matrix established, we can then compute the shortest distance matrix between pairs of doors using the Floyd-Warshall algorithm. The objective function can then be calculated as
\(Z = \sum_{i=1}^N \sum_{j=1}^N w_{ij} d_{ij}\), where \( Z \) is the total weighted door-to-door distance, 
    \( N \) is the total number of rectangles, \( w_{ij} \): flow between door of rectangle \( i \) and door of rectangle \( j \), representing the frequency or cost of travel and \( d_{ij} \) is the shortest path distance between door \( i \) and door \( j \), calculated using the adjacency graph and Euclidean distances.
    
\subsection{Selected Algorithms}

To address the Open Loop Layout Problem (OLLP), we apply four metaheuristic optimization algorithms. These algorithms have been chosen for their ability to effectively navigate the vast and complex solution space of combinatorial problems. Below, we present a brief overview of each algorithm, with more detailed information available in the original studies cited in the respective sections.
\subsubsection{Simple Genetic Algorithm (SGA)}
    The Simple Genetic Algorithm (SGA) is a population-based evolutionary algorithm introduced by \citet{holland1973genetic}. It mimics the process of natural selection through the iterative application of selection (Chooses individuals based on fitness), crossover (Combines parent solutions to generate offspring), and mutation (Introduces random variations to maintain diversity) operators. SGA balances exploration and exploitation to identify high-quality solutions. The pseudocode for the SGA is illustrated in Algorithm \ref{alg:sga}.

\begin{algorithm}
\caption{Simple Genetic Algorithm (SGA)}
\label{alg:sga}
\begin{algorithmic}
\Require Population size $N$, Crossover rate $p_c$, Mutation rate $p_m$, Maximum generations $G_{max}$, Selection method $S$
\Ensure Best solution found
\State Initialize population $P_0$ of size $N$ randomly
\State Set generation counter $g \gets 0$
\While{$g < G_{max}$ and stopping condition not met}
    \State Evaluate fitness of individuals in $P_g$
    \State Select parents from $P_g$ using selection method $S$
    \State Apply crossover with probability $p_c$ to generate offspring
    \State Apply mutation with probability $p_m$ to offspring
    \State Replace old population $P_g$ with new population $P_{g+1}$
    \State Update generation counter $g \gets g + 1$
\EndWhile
\State Return the best solution in the final population
\end{algorithmic}
\end{algorithm}

\subsubsection{Particle Swarm Optimization (PSO)} 
Particle Swarm Optimization, introduced by \citet{kennedy1995pso}, is a swarm intelligence algorithm inspired by the social behavior of birds and fish. Each particle updates its velocity (by combining the particle's inertia, personal best, and global best information) and position (by moving the particle based on its updated velocity) in the solution space based on its own experience and that of its neighbors. A general pseudocode of the PSO is provided in Algorithm~\ref{alg:pso}.
\begin{algorithm}[ht]
\caption{Particle Swarm Optimization (PSO)}
\label{alg:pso}
\begin{algorithmic}[1]  
    \Require Population size $N$, Max generations $G_{max}$, Inertia weight $w$, Cognitive coefficient $c_1$, Social coefficient $c_2$
    \Ensure Best solution found
    \State Initialize positions $\mathbf{x}_i$ and velocities $\mathbf{v}_i$ of $N$ particles randomly
    \State Initialize personal best positions $\mathbf{p}_i$ and global best position $\mathbf{g}$
    \State Set generation counter $g \gets 0$
    \While{$g < G_{max}$ and stopping condition not met}
        \State Evaluate fitness $f(\mathbf{x}_i)$ for each particle $i$
        \State Update personal best positions $\mathbf{p}_i$ based on fitness
        \State Update global best position $\mathbf{g}$ based on the best fitness
        \For{each particle $i$}
            \State Update velocity: 
            $
            \mathbf{v}_i \gets w \mathbf{v}_i + c_1 r_1 (\mathbf{p}_i - \mathbf{x}_i) + c_2 r_2 (\mathbf{g} - \mathbf{x}_i)
            $
            \State Update position: 
            $
            \mathbf{x}_i \gets \mathbf{x}_i + \mathbf{v}_i
            $
        \EndFor
        \State Increment generation counter $g \gets g + 1$
    \EndWhile
    \State \Return the global best solution $\mathbf{g}$
\end{algorithmic}
\end{algorithm}

\subsubsection{Differential Evolution (DE)} 
Differential Evolution, proposed by \citet{storn1997differential}, is a population-based optimization algorithm that employs differential mutation to explore the search space. It is simple yet highly effective for continuous optimization problems. \textbf{Mutation} operator adds the weighted difference of two individuals to a third individual. The
       \textbf{Crossover} operator combines the mutated vector with the target vector. And the \textbf{Selection} operator retains the individual with the best fitness in the next generation. The pseudocode of DE is provided in Algorithm~\ref{alg:DE}.

\begin{algorithm}[ht]
    \caption{Differential Evolution (DE)}
    \label{alg:DE}
    \begin{algorithmic}[1]
        \State \textbf{Input:} Population size $N$, maximum generations $G$, mutation factor $F \in [0, 2]$, crossover rate $CR \in [0, 1]$
        \State \textbf{Output:} Best solution found
        \State Initialize population $\{\mathbf{x}_i^0 \mid i = 1, \dots, N\}$ randomly within the bounds of the search space
        \State Evaluate fitness $f(\mathbf{x}_i^0)$ for all individuals in the initial population
        \For{generation $g = 1$ to $G$}
            \For{each individual $\mathbf{x}_i^g$ in the population}
            \State Select three distinct indices $r_1, r_2, r_3 \neq i$
                \State \textbf{Compute mutant vector} as $ \mathbf{v}_i^g = \mathbf{x}_{r1}^g + F \cdot (\mathbf{x}_{r2}^g - \mathbf{x}_{r3}^g) $

                \State \textbf{Perform crossover}: $ u_{ij}^g = 
    \begin{cases} 
        v_{ij}^g, & \text{if } r_j \leq CR \text{ or } j = j_{\text{rand}} \\
        x_{ij}^g, & \text{otherwise}
    \end{cases}$
                
                \State \textbf{Selection:} Compare trial and target vectors\footnotemark[3]:
                $
                \mathbf{x}_i^{g+1} = 
                \begin{cases} 
                \mathbf{u}_i^g & \text{if } f(\mathbf{u}_i^g) \leq f(\mathbf{x}_i^g) \\
                \mathbf{x}_i^g & \text{otherwise}
                \end{cases}
                $
                
            \EndFor
        \EndFor
        \State \Return the best solution $\mathbf{x}^*$ from the final population
    \end{algorithmic}
\end{algorithm}

\subsubsection{Self-Adaptive Differential Evolution (SADE)}
Self-Adaptive Differential Evolution (SADE), introduced by \citet{OmrSalEng2005SADE}, enhances the classic Differential Evolution by dynamically adjusting mutation and crossover rates. This self-adaptive mechanism improves convergence performance and reduces the need for parameter tuning. The        \textbf{Mutation} operator creates new solutions by adding weighted differences between individuals to a base vector.      
\textbf{Crossover} combines mutated vectors with the current population. And the \textbf{Self-adaptation} technique updates parameters based on success rates of previous iterations. Algorithm~\ref{alg:sade} shows the pseudocode for the DE.

\begin{algorithm}[ht]
    \caption{Self-Adaptive Differential Evolution (SADE)}
    \label{alg:sade}
    \begin{algorithmic}[1] 
        \State \textbf{Input:} Population size $N$, max generations $G$
        \State \textbf{Output:} Best solution found
        \State Initialize population $\{\mathbf{x}_1^0, \mathbf{x}_2^0, \dots, \mathbf{x}_N^0\}$ randomly
        \State Initialize mutation factor $F$ and crossover rate $CR$
        \While{stopping condition not met}
            \For{each individual $i \in \{1, \dots, N\}$}
                \State Select three distinct indices $r_1, r_2, r_3 \neq i$
                \State Generate mutant vector: 
                $ \mathbf{v}_i^g = \mathbf{x}_{r_1}^g + F \cdot (\mathbf{x}_{r_2}^g - \mathbf{x}_{r_3}^g) $
                \State Perform crossover:
                $ u_{ij}^g = 
                    \begin{cases} 
                        v_{ij}^g, & \text{if } r_j \leq CR \text{ or } j = j_{\text{rand}} \\
                        x_{ij}^g, & \text{otherwise}
                \end{cases}$
                \State Evaluate fitness $f(\mathbf{u}_i^g)$
                \If{$f(\mathbf{u}_i^g) \leq f(\mathbf{x}_i^g)$}
                    \State Update population: $\mathbf{x}_i^{g+1} = \mathbf{u}_i^g$
                \Else
                    \State Keep old solution: $\mathbf{x}_i^{g+1} = \mathbf{x}_i^g$
                \EndIf
            \EndFor
            \State Adapt mutation factor $F$ and crossover rate $CR$ based on success rates
        \EndWhile
        \State \Return best solution $\mathbf{x}_{\text{best}}$
    \end{algorithmic}
\end{algorithm}

\subsubsection{Hybrid Methods}
In our pursuit to enhance the efficacy of EAs, we also consider experiments with integrating exact methods into the evolutionary framework, exploring several hybrid approaches:

\begin{itemize}
    \item \textbf{Optimal Placement of Rectangles:} Here, we attempt to optimally position each rectangle individually given a randomly generated entry order. This involved solving a mathematical model aimed at minimizing the weighted distance from already placed rectangles. However, this approach consistently yielded inferior results compared to standalone EAs in our preliminary experiments and thus is not included in the results.
		
    \item \textbf{Individual Rectangle Optimization:} After positioning all rectangles, we investigate the potential improvement of the solution quality by changing the position of individual rectangles. 
\end{itemize}
\noindent	
Due to space constraints, a detailed discussion of hybrid methods is omitted; however, supplementary results are available in the repository at  \url{www.to-be-updated-after-peer-review.com} for interested readers.
\section{Experimental Setup}\label{exp_setup}
\textbf{Benchmark Problems:}
We conducted our experiments on a set of existing benchmark problems with the number of rectangles $\in \{8, 10, 12, 14, 16, 18, 20, 25, 30\}$ acquired from \citep{das, GA}.
Each experiment was assigned a maximum CPU time of $24$ hours, which is a traditional limit in this field and for the benchmark problems used. Each experiment was run $40$ times with different random seed values. All experiments were executed on a HPC server equipped with 16 cores, 
32 threads, and 192 gigabyte of RAM. The experiments are designed to evaluate algorithmic performance, scalability and feasibility.

\noindent
\textbf{Parameter Tuning Selected Algorithms:}
To make a fair comparison and get the best results, we tune the parameters of the selected algorithms to identify the optimal parameters that consistently deliver superior performance across various problem instances. Parameter tuning is a pivotal step in the optimization process, aimed at enhancing the performance of metaheuristic algorithms.

For parameter tuning process we employed the Optuna framework, an efficient hyperparameter optimization tool \citep{akiba2019optuna}. Optuna provides a flexible and lightweight environment for defining and optimizing hyperparameters using a range of search algorithms and pruning strategies. The tuning was performed across a selection of benchmark problems with 8, 10, and 12 rectangles. This approach ensured that the chosen parameters offered robust performance across different problem instances rather than being tailored to a specific instance. The number of generations was capped at $50$ across all algorithms to accelerate the tuning process, with performance evaluated based on the average objective function value across multiple trials. The parameters to be tuned for each algorithm and the best parameter combinations found and used in the experiments are illustrated in the Table~\ref{table:tuned_params} .

 \begin{table}[ht]
 \caption{Optimal parameters for PSO, SADE, SGA, and DE algorithms.}
 \label{table:tuned_params}
\centering
\begin{tabular}{llp{7cm}l}
\hline
\textbf{Algo} & \textbf{Parameter} & \textbf{Definition} & \textbf{Value} \\ \hline
\multirow{6}{*}{PSO}           & $\omega$          & Inertia weight, controlling the influence of the previous velocity          & 0.51  \\ \cmidrule{3-4} 
                                & $\eta_1$          & Cognitive parameter, influencing the particle's tendency to return to its best-known position & 2.42  \\ \cmidrule{3-4} 
                                & $\eta_2$          & Social parameter, influencing the particle's tendency to move towards the swarm's best-known position & 2.37  \\ \cmidrule{3-4} 
                                & Max\_vel          & Maximum velocity, limiting the particle's speed                             & 0.31  \\ \cmidrule{3-4} 
                                & Neighb\_type      & Neighborhood type, determining the topology of particle interactions        & 2     \\ \cmidrule{3-4} 
                                & Neighb\_param     & Neighborhood parameter, specifying the size of the neighborhood             & 4     \\ \hline
\multirow{1}{*}{SADE}           & Variant           & DE strategy variant                                                        & 1     \\ \hline
\multirow{3}{*}{SGA}            & Crossover         & Type of crossover operator                                                  & Binomial  \\ \cmidrule{3-4} 
                                & Mutation          & Type of mutation operator                                                   & Uniform   \\ \cmidrule{3-4} 
                                & Selection         & Method used to select individuals for reproduction                          & Tournament \\ \hline
\multirow{3}{*}{DE}             & F                 & Differential weight, controlling the amplification of differential variations & 0.11   \\ \cmidrule{3-4} 
                                & CR                & Crossover probability, determining the probability of recombination         & 0.86    \\ \cmidrule{3-4} 
                                & Variant           & DE strategy variant                                                        & 9       \\ \hline
                                    \hline
\end{tabular}
\end{table}

\noindent
\textbf{Implementations: }
All the implementations, benchmark problems, and solution methods were implemented in Python 3.7.6. For evolutionary algorithms we used the Pygmo library version 2.16.0~\citep{BisIzzYam2010:pagmo}. 
We used Optuna library in python for parameter tuning \citep{akiba2019optuna}. We also utilized the CPLEX optimization library. 
To encourage further research, we have made our code publicly available at \texttt{www.to-be-updated-after-peer-review.com}.
	
\section{Experimental Results}\label{experimental_results}

The experimental results are benchmarked against the best-known solutions from the literature and are presented in Table \ref{table:results}. For each problem size $n$, the table reports the average objective function value over $40$ runs, with the best solution found across those runs enclosed in curly brackets. The best objective value for each problem instance is highlighted in bold.

\begin{table}
\centering
\caption{Performance comparison of evolutionary algorithms. The format X(Y) represents the Average(Best) objective value. \textbf{N} denotes the number of rectangles.}
\label{table:results}
\begin{tabular}{@{}lccccc@{}}
\toprule
\textbf{N} & \textbf{Literature} & \textbf{DE}                  & \textbf{PSO}                 & \textbf{SADE}                & \textbf{SGA}                 \\
\midrule
8                    & 8906       & 6583 (5617)          & 6665 (5590)          & 6822 (\textbf{5016}) & 6558 (5152)          \\
10                   & 12111      & 13946 (12008)        & 14600 (11937)        & 14585 (\textbf{10712}) & 14343 (10934)        \\
12                   & 36677      & 36434 (31674)        & 36567 (31542)        & 35526 (31233)        & 36302 (\textbf{30239})\\
14                   & 41691      & 44290 (42393)        & 45478 (42131)        & 42855 (40689)        & 39642 (\textbf{36693})\\
16                   & 55064      & 58244 (53763)        & 61733 (55035)        & 57370 (55519)        & 51052 (\textbf{48386})\\
18                   & 66489      & 80770 (76474)        & 81147 (67949)        & 65923 (65179)        & 64451 (\textbf{61330})\\
20                   & -          & 127275 (115777)      & 126114 (108435)      & 117537 (100901)      & 102672 (\textbf{94828})\\
25                   & -          & 278156 (266444)      & 283498 (256110)      & 283330 (257964)      & 234835 (\textbf{223442})\\
30                   & -          & 486497 (469408)      & 494729 (467387)      & 481058 (464196)      & 411587 (\textbf{403351})\\
\bottomrule
\end{tabular}
\end{table}

The results demonstrate that our optimizer consistently delivers competitive outcomes, surpassing the best-known results documented in the literature for problem sizes where comparable data exists. GA proved to be more efficient in solving most problems. This superiority is more prevalent in problems with larger sizes. It is also noteworthy that the optimal objective value for the benchmark problem of size $10$ reported in the literature was achieved through an extensive computational effort involving human-computer collaboration, exceeding a 24-hour time constraint \citep{bv2014}. Despite this, our approach yields superior results.
For larger problem sizes ($n = 20, 25, 30$), where literature results are unavailable, our optimizer effectively achieves feasible solutions within the specified time frame. We invite the research community to utilize the benchmark problems presented in this study as a standard for evaluating and comparing their methods, thereby fostering clarity and collaboration in the field.\\

\begin{figure}[htbp]
    \centering
    \includegraphics[width=0.9\textwidth]{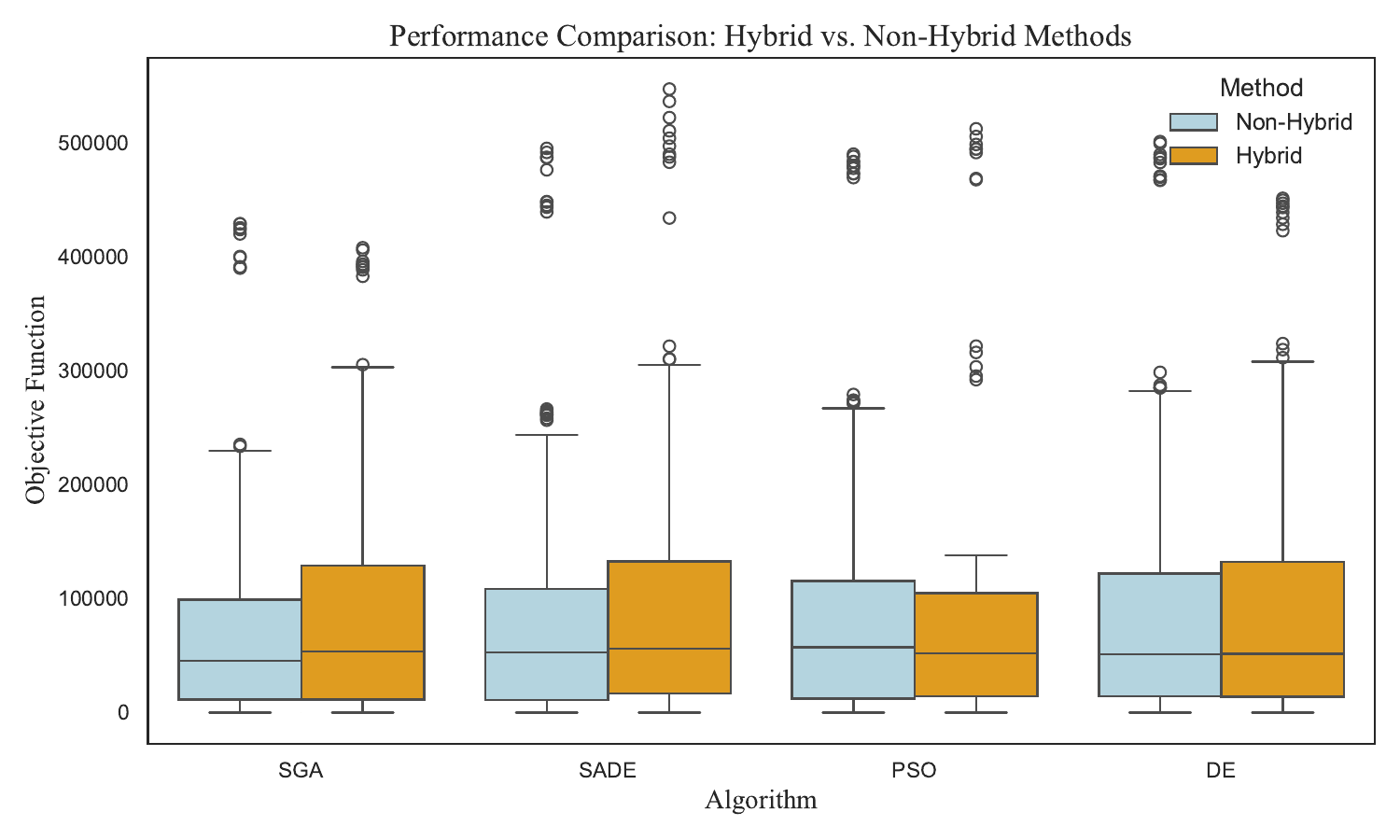}
    \caption{Comparison of algorithm performance with and without hybridization using exact methods.}
    \label{fig:hybrid_nonhybrid_boxplot}
\end{figure}

Figure~\ref{fig:hybrid_nonhybrid_boxplot} illustrates the impact of integrating exact solution methods into evolutionary algorithms. In most cases, hybridization leads to a deterioration in solution quality and increased variance, suggesting that the introduction of exact methods disrupts the natural search dynamics of the algorithms. The only exception is PSO, where the hybrid approach enhances performance, likely due to improved guidance in the particle updates. These observations align with the results in Table~\ref{table:results}, which highlight the superior performance of SGA across various problem sizes. The findings suggest that while hybridization with exact methods may provide benefits in some scenarios and is capable of improving the quality of individuals, it generally introduces instability and does not guarantee improved convergence, particularly for algorithms that rely on iterative evolvement. This is inline with our preliminary observations where using the exact methods for improving the quality of an individual withing initial iterations showed promise, subsequent adjustments failed to produce significant improvements, rendering the results inferior compared to pure EA outcomes.

\begin{figure}[htbp]
    \centering
    \includegraphics[width=0.9\textwidth]{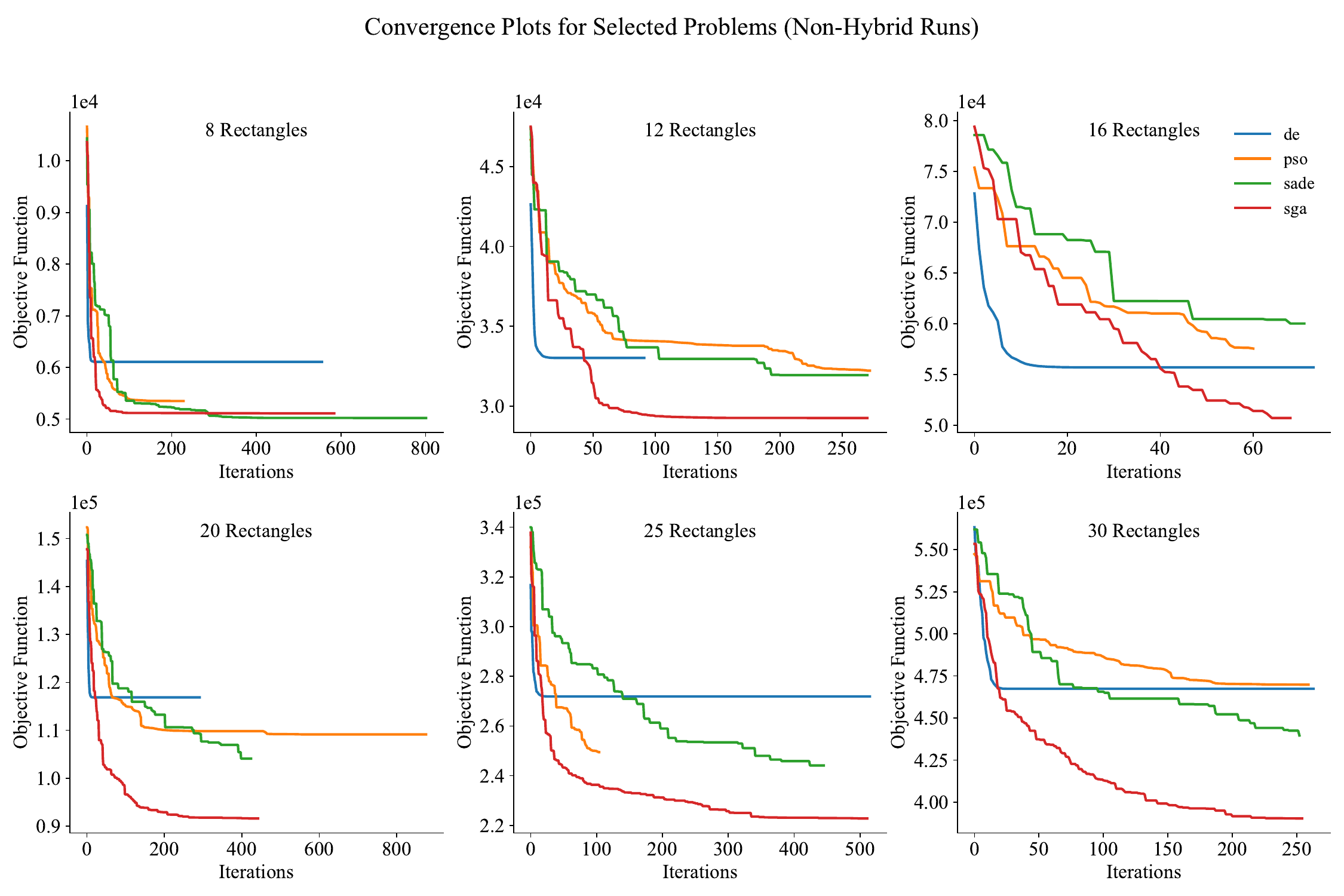}
    \caption{Convergence analysis of evolutionary algorithms for different problem sizes. }
    \label{fig:multi_convergence}
\end{figure}

The convergence analysis illustrated in Figure~\ref{fig:multi_convergence}, highlights distinct performance trends among the evolutionary algorithms across six problem instances (8, 12, 16, 20, 25, and 30 rectangles). DE consistently finds high-quality solutions early but stagnates quickly, suggesting its suitability for time-constrained scenarios, particularly in small-to-medium-sized problems. In contrast, SGA demonstrates a steady improvement over iterations, benefiting from prolonged evolution and population refinement. While the rate of improvement diminishes beyond a certain point, SGA's exploratory capability remains evident. PSO, on the other hand, exhibits fewer total iterations in some cases, indicating a slower internal update mechanism, likely influenced by inertia, velocity updates, and particle interactions. Overall, SGA outperforms the others in terms of both convergence speed and final solution quality, making it a robust choice for layout optimization.

\begin{figure}[htbp]
    \centering
    \includegraphics[width=0.9\textwidth]{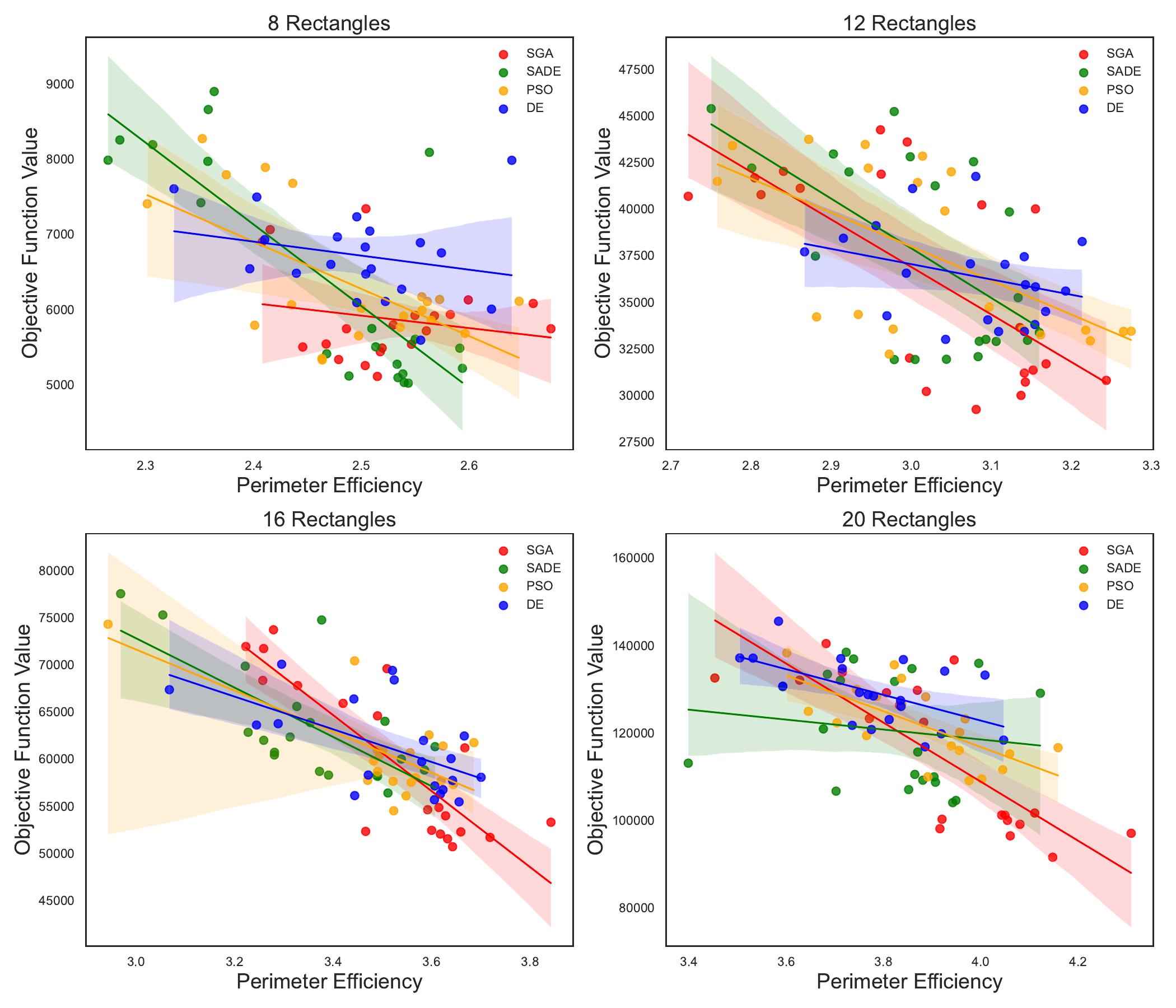}
    \caption{Perimeter efficiency versus objective function. }
    \label{fig:perimeter_efficiency}
\end{figure}

Perimeter efficiency is a critical metric in layout optimization, directly influencing spatial utilization, design effectiveness, and energy performance \citep{MalMueHel1982perimeter}. 
\citet{griffel2020agricultural} define perimeter efficiency as the ratio of the layout perimeter and the occupied area. Figure~\ref{fig:perimeter_efficiency} demonstrates a clear negative correlation between perimeter efficiency and the objective function value across selected problem sizes (8, 12, 16, and 20 rectangles). This indicates that as the objective function improves, perimeter efficiency also increases, emphasizing the trade-off between these metrics. SGA consistently outperforms other algorithms, achieving the highest perimeter efficiency and the lowest objective function values while maintaining robust performance across problem sizes, particularly as complexity increases. In contrast, other algorithms perform competitively in smaller problems but rapidly fall behind as problem size grows. SADE exhibits higher variance, while SGA's results are notably more consistent, with smaller deviations, reinforcing its reliability for both small and large problem instances.

\begin{figure}[htbp]
    \centering
   \includegraphics[width=0.9\textwidth]{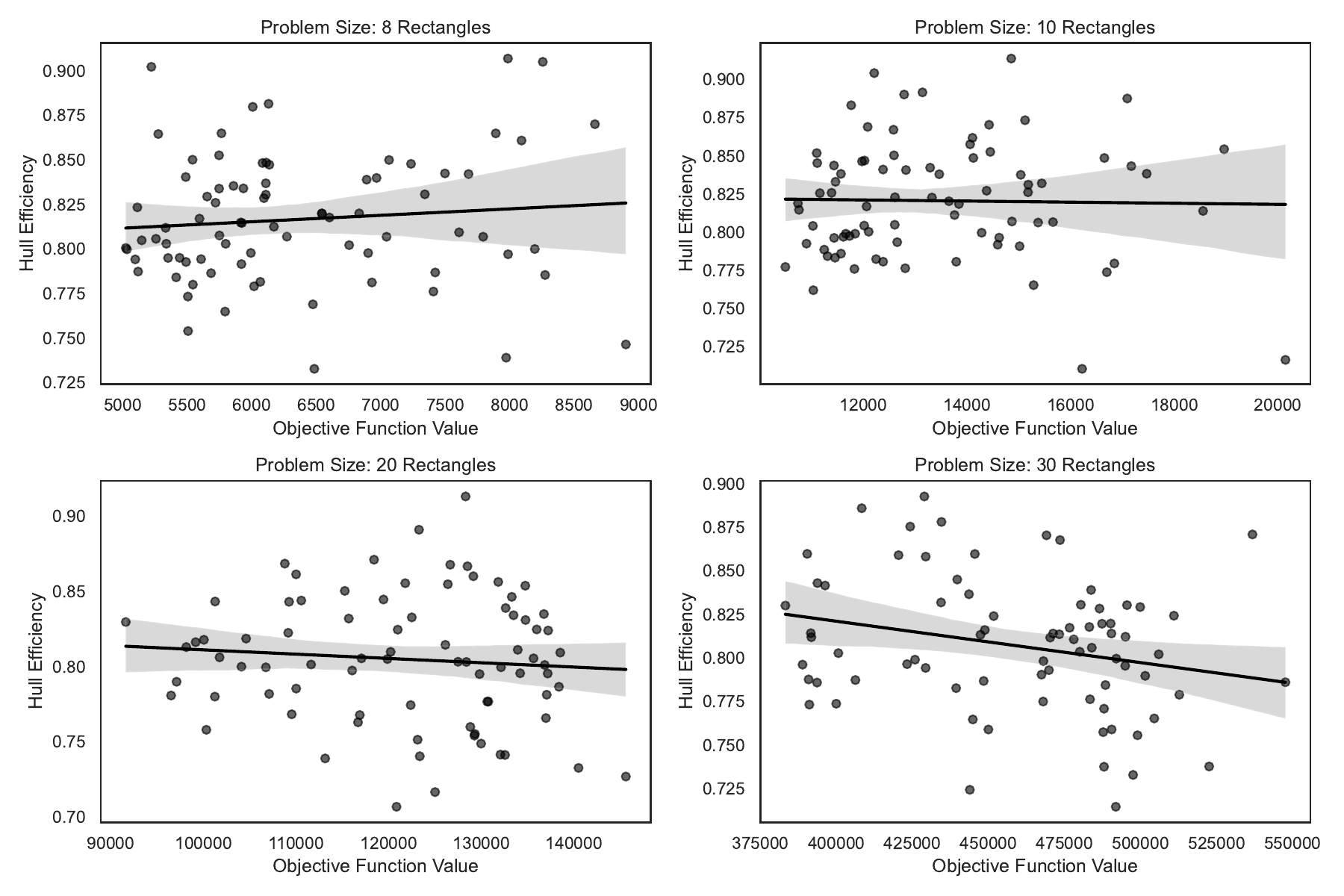}
    \caption{Hull efficiency versus objective function. }
    \label{fig:hull_efficiency_vs_obj}
\end{figure}

Finally, Figure~\ref{fig:hull_efficiency_vs_obj} illustrates the relationship between hull efficiency and the objective function value across different problem sizes (8, 10, 20, and 30 rectangles). Hull efficiency is defined as:
\begin{equation*}
HE = \frac{A_{\text{hull}}}{A_{\text{min bounding box}}}
\end{equation*}
Where: HE is the Hull Efficiency, $A_{\text{hull}}$ is the hull area and $A_{\text{min bounding box}}$ is the area of the smallest rectangle that can accomodate the layout.

Each subplot  in Figure~\ref{fig:hull_efficiency_vs_obj} represents a specific problem size, highlighting how hull efficiency evolves as the objective function improves. A clear downward trend is expected across all cases, indicating that as the objective function value decreases, hull efficiency improves, suggesting a more compact and efficient layout. Interestingly, this is not the case in the problem with $8$ rectangles, where a slight downward trend is observed. This highlights the complexity of open layout problems and the precision required to effectively handle such configurations. The variation in trends among problem sizes also suggests that certain algorithms perform more effectively in maintaining high hull efficiency as problem complexity increases. 

\section{Discussion}\label{discussion}
Evolutionary algorithms (EAs) have proven effective in addressing a broad range of optimization problems, from test suite generation \citep{Dorigo199753} to mechanism synthesis \citep{HUANG2022107928}. Various EAs, such as SGA, DE and PSO, each demonstrate unique strengths in overcoming diverse optimization challenges, and this study reaffirms their applicability to open loop layout problems. We also explored a hierarchical optimization process involving three stages: optimizing the sequence of cells, their orientation, and their positioning. Techniques ranged from sequential exact methods to adjustments guided by EAs. In every instance, pure EAs outperformed these hybrid strategies, underscoring their superior capability in configuring efficient layouts.

The implications of this study are substantial, affecting cost efficiency across multiple domains. For instance, In manufacturing, combining Lean methods with facility layout design can lead to substantial improvements in efficiency, productivity, and cost reduction \citep{Kovács18052020}. Architecturally, optimizing the layout of facilities enhances space utilization and workflow efficiency, applicable in floor planning and other spatial arrangement challenges. Similar benefits extend to logistics through improved packing strategies and to project management via optimized resource allocation. Even in robotic path planning and network design, efficient layouts can minimize operational costs and enhance system performance.

Despite the successes reported here, potential improvements remain unexplored and will be pursued in future studies. Particularly, in large-scale problems, each additional cell can substantially lengthen computation times, especially in calculating the shortest path matrices. Accelerating EA processes through surrogate modeling or alternative techniques such as parallel computation could dramatically enhance computational efficiency and solution quality within fixed time constraints. 

\section{Conclusion}\label{conclusion}
In this paper, we introduced a comprehensive mathematical model for the open loop layout problem, incorporating the novel concept of door-to-door paths. This concept enhances the model's capacity to simulate real-life scenarios more accurately, albeit at the cost of increased computational complexity. We proposed an innovative optimization methodology designed to efficiently arrange and position rectangles within a defined plane, accounting for both geometric and adjacency constraints. The core of our approach lies in a novel encoding scheme that effectively integrates with evolutionary algorithms (EAs).

Extensive experimentation was conducted on a variety of benchmark problems, demonstrating the superiority of our proposed approach over existing methods in the literature. The experimental results indicate that the proposed optimizer consistently yields superior solutions across different problem instances. Additionally, we explored the integration of exact methods at strategic phases of the optimization process. Although these hybrid methods initially exhibited some improvements during the early stages of optimization, their overall performance proved inferior compared to the standalone evolutionary algorithms. Nevertheless, these observations provide valuable insights into the role of exact methods in potentially accelerating initial convergence.
	
Potential applications of our approach span various fields, including VLSI design, architectural floorplanning, packing problems, resource allocation, robotics, and network design. The significant implications of our method include improved efficiency, enhanced performance, cost reduction, scalability, and fostering innovation and competitiveness.
	
Future work will focus on further refining the hybrid optimization approach, exploring adaptive strategies to dynamically balance between heuristic and exact methods based on problem characteristics. One limitation that could be further explored is the rotation of the rectangles, which is currently restricted to only four options. Additionally, extending the approach to three-dimensional packing and exploring other practical applications will be valuable directions for further research. Finally, application of machine learning techniques to implement surrogate models for improving the efficiency of the proposed method in large-scale problems will be interesting. To facilitate the contributions in the field we make our algorithms and data publicly available at \url{www.to-be-updated-after-peer-review.com}.

\section*{Statements and Declarations}

\textbf{Funding} \\
No funds, grants, or other support was received.

\noindent
\textbf{Conflict of Interest} \\
The authors declare that they have no conflicts of interest. The authors have no relevant financial or non-financial interests to disclose.

\noindent
\textbf{Ethics Approval and Consent to Participate} \\
Not applicable.\\
\noindent
\textbf{Consent for Publication} \\
All authors have given their consent for the publication of this study.\\
\noindent
\textbf{Code and Data Availability} \\
The implemented code and datasets generated and analyzed during this study will be made available after peer review. To comply with the journal’s double-blind review process, links to these resources have been temporarily replaced with placeholders.
\noindent
\textbf{Materials Availability} \\
Not applicable.
\noindent

\bibliography{main.bbl}


\begin{thebibliography}{53}
\ifx \bisbn   \undefined \def \bisbn  #1{ISBN #1}\fi
\ifx \binits  \undefined \def \binits#1{#1}\fi
\ifx \bauthor  \undefined \def \bauthor#1{#1}\fi
\ifx \batitle  \undefined \def \batitle#1{#1}\fi
\ifx \bjtitle  \undefined \def \bjtitle#1{#1}\fi
\ifx \bvolume  \undefined \def \bvolume#1{\textbf{#1}}\fi
\ifx \byear  \undefined \def \byear#1{#1}\fi
\ifx \bissue  \undefined \def \bissue#1{#1}\fi
\ifx \bfpage  \undefined \def \bfpage#1{#1}\fi
\ifx \blpage  \undefined \def \blpage #1{#1}\fi
\ifx \burl  \undefined \def \burl#1{\textsf{#1}}\fi
\ifx \doiurl  \undefined \def \doiurl#1{\url{https://doi.org/#1}}\fi
\ifx \betal  \undefined \def \betal{\textit{et al.}}\fi
\ifx \binstitute  \undefined \def \binstitute#1{#1}\fi
\ifx \binstitutionaled  \undefined \def \binstitutionaled#1{#1}\fi
\ifx \bctitle  \undefined \def \bctitle#1{#1}\fi
\ifx \beditor  \undefined \def \beditor#1{#1}\fi
\ifx \bpublisher  \undefined \def \bpublisher#1{#1}\fi
\ifx \bbtitle  \undefined \def \bbtitle#1{#1}\fi
\ifx \bedition  \undefined \def \bedition#1{#1}\fi
\ifx \bseriesno  \undefined \def \bseriesno#1{#1}\fi
\ifx \blocation  \undefined \def \blocation#1{#1}\fi
\ifx \bsertitle  \undefined \def \bsertitle#1{#1}\fi
\ifx \bsnm \undefined \def \bsnm#1{#1}\fi
\ifx \bsuffix \undefined \def \bsuffix#1{#1}\fi
\ifx \bparticle \undefined \def \bparticle#1{#1}\fi
\ifx \barticle \undefined \def \barticle#1{#1}\fi
\bibcommenthead
\ifx \bconfdate \undefined \def \bconfdate #1{#1}\fi
\ifx \botherref \undefined \def \botherref #1{#1}\fi
\ifx \url \undefined \def \url#1{\textsf{#1}}\fi
\ifx \bchapter \undefined \def \bchapter#1{#1}\fi
\ifx \bbook \undefined \def \bbook#1{#1}\fi
\ifx \bcomment \undefined \def \bcomment#1{#1}\fi
\ifx \oauthor \undefined \def \oauthor#1{#1}\fi
\ifx \citeauthoryear \undefined \def \citeauthoryear#1{#1}\fi
\ifx \endbibitem  \undefined \def \endbibitem {}\fi
\ifx \bconflocation  \undefined \def \bconflocation#1{#1}\fi
\ifx \arxivurl  \undefined \def \arxivurl#1{\textsf{#1}}\fi
\csname PreBibitemsHook\endcsname

\bibitem[\protect\citeauthoryear{Alduaij and Hassan}{2020}]{alduaij2020adopting}
\begin{barticle}
\bauthor{\bsnm{Alduaij}, \binits{A.}},
\bauthor{\bsnm{Hassan}, \binits{N.M.}}:
\batitle{Adopting a circular open-field layout in designing flexible manufacturing systems}.
\bjtitle{International Journal of Computer Integrated Manufacturing}
\bvolume{33}(\bissue{6}),
\bfpage{572}--\blpage{589}
(\byear{2020})
\end{barticle}
\endbibitem

\bibitem[\protect\citeauthoryear{Akiba et~al.}{2019}]{akiba2019optuna}
\begin{botherref}
\oauthor{\bsnm{Akiba}, \binits{T.}},
\oauthor{\bsnm{Sano}, \binits{S.}},
\oauthor{\bsnm{Yanase}, \binits{T.}},
\oauthor{\bsnm{Ohta}, \binits{T.}},
\oauthor{\bsnm{Koyama}, \binits{M.}}:
Optuna: A next-generation hyperparameter optimization framework.
Proceedings of the 25th ACM SIGKDD International Conference on Knowledge Discovery \& Data Mining,
2623--2631
(2019).
ACM
\end{botherref}
\endbibitem

\bibitem[\protect\citeauthoryear{Benson and Foote}{1997}]{benson1997doorfast}
\begin{barticle}
\bauthor{\bsnm{Benson}, \binits{B.}},
\bauthor{\bsnm{Foote}, \binits{B.L.}}:
\batitle{Doorfast: A constructive procedure to optimally layout a facility including aisles and door locations basedon anaisle flow distance metric}.
\bjtitle{International Journal of Production Research}
\bvolume{35}(\bissue{7}),
\bfpage{1825}--\blpage{1842}
(\byear{1997})
\end{barticle}
\endbibitem

\bibitem[\protect\citeauthoryear{Biscani et~al.}{2010}]{BisIzzYam2010:pagmo}
\begin{bchapter}
\bauthor{\bsnm{Biscani}, \binits{F.}},
\bauthor{\bsnm{Izzo}, \binits{D.}},
\bauthor{\bsnm{Yam}, \binits{C.H.}}:
\bctitle{A global optimisation toolbox for massively parallel engineering optimisation}.
(\byear{2010}).
\burl{http://arxiv.org/abs/1004.3824}
\end{bchapter}
\endbibitem

\bibitem[\protect\citeauthoryear{Bukchin and Tzur}{2014}]{bukTzu2014new}
\begin{barticle}
\bauthor{\bsnm{Bukchin}, \binits{Y.}},
\bauthor{\bsnm{Tzur}, \binits{M.}}:
\batitle{A new milp approach for the facility process-layout design problem with rectangular and l/t shape departments}.
\bjtitle{International Journal of Production Research}
\bvolume{52}(\bissue{24}),
\bfpage{7339}--\blpage{7359}
(\byear{2014})
\end{barticle}
\endbibitem

\bibitem[\protect\citeauthoryear{Carlo and Bozer}{2011}]{carlo2011analysis}
\begin{barticle}
\bauthor{\bsnm{Carlo}, \binits{H.J.}},
\bauthor{\bsnm{Bozer}, \binits{Y.A.}}:
\batitle{Analysis of optimum shape and door assignment problems in rectangular unit-load crossdocks}.
\bjtitle{International Journal of Logistics Research and Applications}
\bvolume{14}(\bissue{3}),
\bfpage{149}--\blpage{163}
(\byear{2011})
\end{barticle}
\endbibitem

\bibitem[\protect\citeauthoryear{Chae and Peters}{2006}]{Chae01072006}
\begin{barticle}
\bauthor{\bsnm{Chae}, \binits{J.}},
\bauthor{\bsnm{Peters}, \binits{B.A.}}:
\batitle{A simulated annealing algorithm based on a closed loop layout for facility layout design in flexible manufacturing systems}.
\bjtitle{International Journal of Production Research}
\bvolume{44}(\bissue{13}),
\bfpage{2561}--\blpage{2572}
(\byear{2006})
\doiurl{10.1080/00207540500446287}
{\href{https://arxiv.org/abs/https://doi.org/10.1080/00207540500446287}{{https://doi.org/10.1080/00207540500446287}}}
\end{barticle}
\endbibitem

\bibitem[\protect\citeauthoryear{Das}{1993}]{das}
\begin{barticle}
\bauthor{\bsnm{Das}, \binits{S.K.}}:
\batitle{A facility layout method for feasible manufacturing systems}.
\bjtitle{International Journal of Production Research}
\bvolume{31}(\bissue{2}),
\bfpage{279}--\blpage{297}
(\byear{1993})
\end{barticle}
\endbibitem

\bibitem[\protect\citeauthoryear{Dorigo and Gambardella}{1997}]{Dorigo199753}
\begin{barticle}
\bauthor{\bsnm{Dorigo}, \binits{M.}},
\bauthor{\bsnm{Gambardella}, \binits{L.M.}}:
\batitle{Ant colony system: A cooperative learning approach to the traveling salesman problem}.
\bjtitle{IEEE Transactions on Evolutionary Computation}
\bvolume{1}(\bissue{1}),
\bfpage{53}--\blpage{66}
(\byear{1997})
\doiurl{10.1109/4235.585892} .
\bcomment{Cited by: 6935}
\end{barticle}
\endbibitem

\bibitem[\protect\citeauthoryear{De~Guzman et~al.}{1997}]{de1997complexity}
\begin{barticle}
\bauthor{\bsnm{De~Guzman}, \binits{M.}},
\bauthor{\bsnm{Prabhu}, \binits{N.}},
\bauthor{\bsnm{Tanchoco}, \binits{J.}}:
\batitle{Complexity of the agv shortest path and single-loop guide path layout problems}.
\bjtitle{International Journal of Production Research}
\bvolume{35}(\bissue{8}),
\bfpage{2083}--\blpage{2092}
(\byear{1997})
\end{barticle}
\endbibitem

\bibitem[\protect\citeauthoryear{Dorrah and Marzouk}{2021}]{dorrah2021integrated}
\begin{barticle}
\bauthor{\bsnm{Dorrah}, \binits{D.H.}},
\bauthor{\bsnm{Marzouk}, \binits{M.}}:
\batitle{Integrated multi-objective optimization and agent-based building occupancy modeling for space layout planning}.
\bjtitle{Journal of Building Engineering}
\bvolume{34},
\bfpage{101902}
(\byear{2021})
\end{barticle}
\endbibitem

\bibitem[\protect\citeauthoryear{Drira et~al.}{2006}]{drira2006survey}
\begin{barticle}
\bauthor{\bsnm{Drira}, \binits{A.}},
\bauthor{\bsnm{Pierreval}, \binits{H.}},
\bauthor{\bsnm{Hajri-Gabouj}, \binits{S.}}:
\batitle{Facility layout problems: A survey}.
\bjtitle{Annual Reviews in Control}
\bvolume{31}(\bissue{2}),
\bfpage{255}--\blpage{267}
(\byear{2006})
\end{barticle}
\endbibitem

\bibitem[\protect\citeauthoryear{D{\'\i}az et~al.}{2024}]{diaz2024variable}
\begin{botherref}
\oauthor{\bsnm{D{\'\i}az}, \binits{J.}},
\oauthor{\bsnm{Rodriguez}, \binits{H.}},
\oauthor{\bsnm{Fajardo-Calder{\'\i}n}, \binits{J.}},
\oauthor{\bsnm{Angulo}, \binits{I.}},
\oauthor{\bsnm{Onieva}, \binits{E.}}:
A variable neighbourhood search for minimization of operation times through warehouse layout optimization.
Logic Journal of the IGPL,
018
(2024)
\end{botherref}
\endbibitem

\bibitem[\protect\citeauthoryear{Gu et~al.}{2010}]{GU2010539}
\begin{barticle}
\bauthor{\bsnm{Gu}, \binits{J.}},
\bauthor{\bsnm{Goetschalckx}, \binits{M.}},
\bauthor{\bsnm{McGinnis}, \binits{L.F.}}:
\batitle{Research on warehouse design and performance evaluation: A comprehensive review}.
\bjtitle{European Journal of Operational Research}
\bvolume{203}(\bissue{3}),
\bfpage{539}--\blpage{549}
(\byear{2010})
\doiurl{10.1016/j.ejor.2009.07.031}
\end{barticle}
\endbibitem

\bibitem[\protect\citeauthoryear{Griffel et~al.}{2020}]{griffel2020agricultural}
\begin{barticle}
\bauthor{\bsnm{Griffel}, \binits{L.M.}},
\bauthor{\bsnm{Vazhnik}, \binits{V.}},
\bauthor{\bsnm{Hartley}, \binits{D.S.}},
\bauthor{\bsnm{Hansen}, \binits{J.K.}},
\bauthor{\bsnm{Roni}, \binits{M.}}:
\batitle{Agricultural field shape descriptors as predictors of field efficiency for perennial grass harvesting: An empirical proof}.
\bjtitle{Computers and Electronics in Agriculture}
\bvolume{168},
\bfpage{105088}
(\byear{2020})
\end{barticle}
\endbibitem

\bibitem[\protect\citeauthoryear{Hauser and Chung}{2006}]{Hauser01112006}
\begin{barticle}
\bauthor{\bsnm{Hauser}, \binits{K.}},
\bauthor{\bsnm{Chung}, \binits{C.H.}}:
\batitle{Genetic algorithms for layout optimization in crossdocking operations of a manufacturing plant}.
\bjtitle{International Journal of Production Research}
\bvolume{44}(\bissue{21}),
\bfpage{4663}--\blpage{4680}
(\byear{2006})
\doiurl{10.1080/00207540500521147}
\end{barticle}
\endbibitem

\bibitem[\protect\citeauthoryear{Heragu and Kusiak}{1991}]{hk}
\begin{barticle}
\bauthor{\bsnm{Heragu}, \binits{S.S.}},
\bauthor{\bsnm{Kusiak}, \binits{A.}}:
\batitle{Efficient models for facility layout problem}.
\bjtitle{European Journal of Operations Research}
\bvolume{53},
\bfpage{1}--\blpage{13}
(\byear{1991})
\end{barticle}
\endbibitem

\bibitem[\protect\citeauthoryear{Holland}{1973}]{holland1973genetic}
\begin{barticle}
\bauthor{\bsnm{Holland}, \binits{J.H.}}:
\batitle{Genetic algorithms and the optimal allocation of trials}.
\bjtitle{SIAM journal on computing}
\bvolume{2}(\bissue{2}),
\bfpage{88}--\blpage{105}
(\byear{1973})
\doiurl{10.1137/0202009}
\end{barticle}
\endbibitem

\bibitem[\protect\citeauthoryear{Huang et~al.}{2022}]{HUANG2022107928}
\begin{barticle}
\bauthor{\bsnm{Huang}, \binits{Q.}},
\bauthor{\bsnm{Yu}, \binits{Y.}},
\bauthor{\bsnm{Zhang}, \binits{K.}},
\bauthor{\bsnm{Li}, \binits{S.}},
\bauthor{\bsnm{Lu}, \binits{H.}},
\bauthor{\bsnm{Li}, \binits{J.}},
\bauthor{\bsnm{Zhang}, \binits{A.}},
\bauthor{\bsnm{Mei}, \binits{T.}}:
\batitle{Optimal synthesis of mechanisms using repellency evolutionary algorithm}.
\bjtitle{Knowledge-Based Systems}
\bvolume{239},
\bfpage{107928}
(\byear{2022})
\doiurl{10.1016/j.knosys.2021.107928}
\end{barticle}
\endbibitem

\bibitem[\protect\citeauthoryear{J~A~Bennell and Whitehead}{2002}]{Bennell01102002}
\begin{barticle}
\bauthor{\bsnm{J~A~Bennell}, \binits{C.N.P.}},
\bauthor{\bsnm{Whitehead}, \binits{J.D.}}:
\batitle{Local search algorithms for the min-max loop layout problem}.
\bjtitle{Journal of the Operational Research Society}
\bvolume{53}(\bissue{10}),
\bfpage{1109}--\blpage{1117}
(\byear{2002})
\doiurl{10.1057/palgrave.jors.2601269}
\end{barticle}
\endbibitem

\bibitem[\protect\citeauthoryear{Kundu and Dan}{2012}]{kundu2012metaheuristic}
\begin{barticle}
\bauthor{\bsnm{Kundu}, \binits{A.}},
\bauthor{\bsnm{Dan}, \binits{P.K.}}:
\batitle{Metaheuristic in facility layout problems: current trend and future direction}.
\bjtitle{International Journal of Industrial and Systems Engineering}
\bvolume{10}(\bissue{2}),
\bfpage{238}--\blpage{253}
(\byear{2012})
\end{barticle}
\endbibitem

\bibitem[\protect\citeauthoryear{Kennedy and Eberhart}{1995}]{kennedy1995pso}
\begin{bchapter}
\bauthor{\bsnm{Kennedy}, \binits{J.}},
\bauthor{\bsnm{Eberhart}, \binits{R.}}:
\bctitle{Particle swarm optimization}.
In: \bbtitle{Proceedings of ICNN'95-international Conference on Neural Networks},
vol. \bseriesno{4},
pp. \bfpage{1942}--\blpage{1948}
(\byear{1995}).
\bcomment{ieee}
\end{bchapter}
\endbibitem

\bibitem[\protect\citeauthoryear{Kim and Kim}{2015}]{kim2015optimal}
\begin{barticle}
\bauthor{\bsnm{Kim}, \binits{Y.-D.}},
\bauthor{\bsnm{Kim}, \binits{H.-S.}}:
\batitle{Optimal warehouse design and layout with consideration of material handling system}.
\bjtitle{Computers \& Industrial Engineering}
\bvolume{87},
\bfpage{464}--\blpage{476}
(\byear{2015})
\end{barticle}
\endbibitem

\bibitem[\protect\citeauthoryear{Kang et~al.}{2018}]{kang2018closed}
\begin{barticle}
\bauthor{\bsnm{Kang}, \binits{S.}},
\bauthor{\bsnm{Kim}, \binits{M.}},
\bauthor{\bsnm{Chae}, \binits{J.}}:
\batitle{A closed loop based facility layout design using a cuckoo search algorithm}.
\bjtitle{Expert Systems with Applications}
\bvolume{93},
\bfpage{322}--\blpage{335}
(\byear{2018})
\end{barticle}
\endbibitem

\bibitem[\protect\citeauthoryear{Klausnitzer and Lasch}{2019}]{KlaRai2019OptimalFL}
\begin{barticle}
\bauthor{\bsnm{Klausnitzer}, \binits{A.}},
\bauthor{\bsnm{Lasch}, \binits{R.}}:
\batitle{Optimal facility layout and material handling network design}.
\bjtitle{Comput. Oper. Res.}
\bvolume{103},
\bfpage{237}--\blpage{251}
(\byear{2019})
\end{barticle}
\endbibitem

\bibitem[\protect\citeauthoryear{Kovács}{2020}]{Kovács18052020}
\begin{barticle}
\bauthor{\bsnm{Kovács}, \binits{G.}}:
\batitle{Combination of lean value-oriented conception and facility layout design for even more significant efficiency improvement and cost reduction}.
\bjtitle{International Journal of Production Research}
\bvolume{58}(\bissue{10}),
\bfpage{2916}--\blpage{2936}
(\byear{2020})
\doiurl{10.1080/00207543.2020.1712490}
{\href{https://arxiv.org/abs/https://doi.org/10.1080/00207543.2020.1712490}{{https://doi.org/10.1080/00207543.2020.1712490}}}
\end{barticle}
\endbibitem

\bibitem[\protect\citeauthoryear{Khalilabadi et~al.}{2024}]{khalilabadi2024exploiting}
\begin{botherref}
\oauthor{\bsnm{Khalilabadi}, \binits{S.M.G.}},
\oauthor{\bsnm{Roy}, \binits{D.}},
\oauthor{\bsnm{Koster}, \binits{R.}}:
Exploiting travel sequences to optimise facility layouts with multiple input/output points.
International Journal of Production Research,
1--29
(2024)
\end{botherref}
\endbibitem

\bibitem[\protect\citeauthoryear{Meller and Bozer}{1996}]{MelBoz1996PR}
\begin{barticle}
\bauthor{\bsnm{Meller}, \binits{R.D.}},
\bauthor{\bsnm{Bozer}, \binits{Y.A.}}:
\batitle{A new simulated annealing algorithm for the facility layout problem}.
\bjtitle{International Journal of Production Research}
\bvolume{34}(\bissue{6}),
\bfpage{1675}--\blpage{1692}
(\byear{1996})
\doiurl{10.1080/00207549608904990}
{\href{https://arxiv.org/abs/https://doi.org/10.1080/00207549608904990}{{https://doi.org/10.1080/00207549608904990}}}
\end{barticle}
\endbibitem

\bibitem[\protect\citeauthoryear{Meng et~al.}{2021}]{meng2021impact}
\begin{barticle}
\bauthor{\bsnm{Meng}, \binits{L.}},
\bauthor{\bsnm{Batt}, \binits{R.J.}},
\bauthor{\bsnm{Terwiesch}, \binits{C.}}:
\batitle{The impact of facility layout on service worker behavior: An empirical study of nurses in the emergency department}.
\bjtitle{Manufacturing \& Service Operations Management}
\bvolume{23}(\bissue{4}),
\bfpage{819}--\blpage{834}
(\byear{2021})
\end{barticle}
\endbibitem

\bibitem[\protect\citeauthoryear{Maling et~al.}{1982}]{MalMueHel1982perimeter}
\begin{bchapter}
\bauthor{\bsnm{Maling}, \binits{K.}},
\bauthor{\bsnm{Mueller}, \binits{S.H.}},
\bauthor{\bsnm{Heller}, \binits{W.R.}}:
\bctitle{On finding most optimal rectangular package plans}.
In: \bbtitle{Proceedings of the 19th Design Automation Conference}.
\bsertitle{DAC '82},
pp. \bfpage{663}--\blpage{670}.
\bpublisher{IEEE Press}, \blocation{???}
(\byear{1982})
\end{bchapter}
\endbibitem

\bibitem[\protect\citeauthoryear{Meller et~al.}{1999}]{meller}
\begin{barticle}
\bauthor{\bsnm{Meller}, \binits{R.D.}},
\bauthor{\bsnm{Narazanan}, \binits{V.}},
\bauthor{\bsnm{Vance}, \binits{P.H.}}:
\batitle{Optimal facility layout design}.
\bjtitle{Operations Research Letters}
\bvolume{23},
\bfpage{117}--\blpage{127}
(\byear{1999})
\end{barticle}
\endbibitem

\bibitem[\protect\citeauthoryear{Montreuil}{1990}]{mont}
\begin{botherref}
\oauthor{\bsnm{Montreuil}, \binits{B.}}:
A modeling framework for integrating layout design and flow network design.
Proceedings of the Material Handling Research Colloquium,
43--58
(1990)
\end{botherref}
\endbibitem

\bibitem[\protect\citeauthoryear{Ma et~al.}{2025}]{ma2025enhancing}
\begin{barticle}
\bauthor{\bsnm{Ma}, \binits{C.}},
\bauthor{\bsnm{Song}, \binits{M.}},
\bauthor{\bsnm{Zeng}, \binits{W.}},
\bauthor{\bsnm{Wang}, \binits{X.}},
\bauthor{\bsnm{Chen}, \binits{T.}},
\bauthor{\bsnm{Wu}, \binits{S.}}:
\batitle{Enhancing urban emergency response: A euclidean distance-based framework for optimizing rescue facility layouts}.
\bjtitle{Sustainable Cities and Society}
\bvolume{118},
\bfpage{106006}
(\byear{2025})
\end{barticle}
\endbibitem

\bibitem[\protect\citeauthoryear{Muther}{1973}]{muther1973systematic}
\begin{bbook}
\bauthor{\bsnm{Muther}, \binits{R.}}:
\bbtitle{Systematic Layout Planning}.
\bpublisher{Industrial Press, Inc.},
\blocation{Chicago, IL}
(\byear{1973})
\end{bbook}
\endbibitem

\bibitem[\protect\citeauthoryear{Niroomand and Vizv{\'a}ri}{2013}]{niroomand2013mixed}
\begin{barticle}
\bauthor{\bsnm{Niroomand}, \binits{S.}},
\bauthor{\bsnm{Vizv{\'a}ri}, \binits{B.}}:
\batitle{A mixed integer linear programming formulation of closed loop layout with exact distances}.
\bjtitle{Journal of Industrial and Production Engineering}
\bvolume{30}(\bissue{3}),
\bfpage{190}--\blpage{201}
(\byear{2013})
\end{barticle}
\endbibitem

\bibitem[\protect\citeauthoryear{Omran et~al.}{2005}]{OmrSalEng2005SADE}
\begin{bchapter}
\bauthor{\bsnm{Omran}, \binits{M.G.H.}},
\bauthor{\bsnm{Salman}, \binits{A.}},
\bauthor{\bsnm{Engelbrecht}, \binits{A.P.}}:
\bctitle{Self-adaptive differential evolution}.
In: \beditor{\bsnm{Hao}, \binits{Y.}},
\beditor{\bsnm{Liu}, \binits{J.}},
\beditor{\bsnm{Wang}, \binits{Y.}},
\beditor{\bsnm{Cheung}, \binits{Y.-m.}},
\beditor{\bsnm{Yin}, \binits{H.}},
\beditor{\bsnm{Jiao}, \binits{L.}},
\beditor{\bsnm{Ma}, \binits{J.}},
\beditor{\bsnm{Jiao}, \binits{Y.-C.}} (eds.)
\bbtitle{Computational Intelligence and Security},
pp. \bfpage{192}--\blpage{199}.
\bpublisher{Springer},
\blocation{Berlin, Heidelberg}
(\byear{2005})
\end{bchapter}
\endbibitem

\bibitem[\protect\citeauthoryear{P{\'e}rez-Gosende et~al.}{2021}]{perez2021facility}
\begin{barticle}
\bauthor{\bsnm{P{\'e}rez-Gosende}, \binits{P.}},
\bauthor{\bsnm{Mula}, \binits{J.}},
\bauthor{\bsnm{D{\'\i}az-Madro{\~n}ero}, \binits{M.}}:
\batitle{Facility layout planning. an extended literature review}.
\bjtitle{International Journal of Production Research}
\bvolume{59}(\bissue{12}),
\bfpage{3777}--\blpage{3816}
(\byear{2021})
\end{barticle}
\endbibitem

\bibitem[\protect\citeauthoryear{Powell}{2010}]{powell2010merging}
\begin{barticle}
\bauthor{\bsnm{Powell}, \binits{W.B.}}:
\batitle{Merging ai and or to solve high-dimensional stochastic optimization problems using approximate dynamic programming}.
\bjtitle{INFORMS Journal on Computing}
\bvolume{22}(\bissue{1}),
\bfpage{2}--\blpage{17}
(\byear{2010})
\end{barticle}
\endbibitem

\bibitem[\protect\citeauthoryear{Pourvaziri et~al.}{2021}]{pouPieMar2021integrating}
\begin{barticle}
\bauthor{\bsnm{Pourvaziri}, \binits{H.}},
\bauthor{\bsnm{Pierreval}, \binits{H.}},
\bauthor{\bsnm{Marian}, \binits{H.}}:
\batitle{Integrating facility layout design and aisle structure in manufacturing systems: Formulation and exact solution}.
\bjtitle{European journal of Operational Research}
\bvolume{290}(\bissue{2}),
\bfpage{499}--\blpage{513}
(\byear{2021})
\end{barticle}
\endbibitem

\bibitem[\protect\citeauthoryear{Preparata and Shamos}{1985}]{preparata1985computational}
\begin{bbook}
\bauthor{\bsnm{Preparata}, \binits{F.P.}},
\bauthor{\bsnm{Shamos}, \binits{M.I.}}:
\bbtitle{Computational Geometry: An Introduction}.
\bpublisher{Springer},
\blocation{New York}
(\byear{1985})
\end{bbook}
\endbibitem

\bibitem[\protect\citeauthoryear{Rajasekharan et~al.}{1998}]{GA}
\begin{barticle}
\bauthor{\bsnm{Rajasekharan}, \binits{M.}},
\bauthor{\bsnm{Peters}, \binits{B.A.}},
\bauthor{\bsnm{Yang}, \binits{T.}}:
\batitle{A genetic algorithm for facility layout design in flexible manufacturing systems}.
\bjtitle{Int. J. Prod. Res.}
\bvolume{36},
\bfpage{95}--\blpage{110}
(\byear{1998})
\end{barticle}
\endbibitem

\bibitem[\protect\citeauthoryear{Sherali et~al.}{2003}]{sherali}
\begin{barticle}
\bauthor{\bsnm{Sherali}, \binits{H.D.}},
\bauthor{\bsnm{Fraticelli}, \binits{B.M.P.}},
\bauthor{\bsnm{Meller}, \binits{R.D.}}:
\batitle{Enhanced model formulation for optimal layout}.
\bjtitle{Operations Research}
\bvolume{51},
\bfpage{629}--\blpage{644}
(\byear{2003})
\end{barticle}
\endbibitem

\bibitem[\protect\citeauthoryear{Shi et~al.}{2022}]{shi2022study}
\begin{barticle}
\bauthor{\bsnm{Shi}, \binits{Z.}},
\bauthor{\bsnm{Li}, \binits{Y.}},
\bauthor{\bsnm{Boh{\'a}cs}, \binits{G.}},
\bauthor{\bsnm{Zhou}, \binits{Q.}}:
\batitle{A study on optimal location selection and semi-finished product inventory allocation in the steel industry}.
\bjtitle{Sustainability}
\bvolume{14}(\bissue{22}),
\bfpage{15279}
(\byear{2022})
\end{barticle}
\endbibitem

\bibitem[\protect\citeauthoryear{Storn and Price}{1997}]{storn1997differential}
\begin{barticle}
\bauthor{\bsnm{Storn}, \binits{R.}},
\bauthor{\bsnm{Price}, \binits{K.}}:
\batitle{Differential evolution--a simple and efficient heuristic for global optimization over continuous spaces}.
\bjtitle{Journal of Global Optimization}
\bvolume{11}(\bissue{4}),
\bfpage{341}--\blpage{359}
(\byear{1997})
\end{barticle}
\endbibitem

\bibitem[\protect\citeauthoryear{Thornton et~al.}{1979}]{thornton1979rectangular}
\begin{barticle}
\bauthor{\bsnm{Thornton}, \binits{V.D.}},
\bauthor{\bsnm{Francis}, \binits{R.L.}},
\bauthor{\bsnm{Lowe}, \binits{T.J.}}:
\batitle{Rectangular layout problems with worst-case distance measures}.
\bjtitle{AIIE Transactions}
\bvolume{11}(\bissue{1}),
\bfpage{2}--\blpage{11}
(\byear{1979})
\end{barticle}
\endbibitem

\bibitem[\protect\citeauthoryear{Tongur et~al.}{2020}]{TONGUR2020951}
\begin{barticle}
\bauthor{\bsnm{Tongur}, \binits{V.}},
\bauthor{\bsnm{Hacibeyoglu}, \binits{M.}},
\bauthor{\bsnm{Ulker}, \binits{E.}}:
\batitle{Solving a big-scaled hospital facility layout problem with meta-heuristics algorithms}.
\bjtitle{Engineering Science and Technology, an International Journal}
\bvolume{23}(\bissue{4}),
\bfpage{951}--\blpage{959}
(\byear{2020})
\doiurl{10.1016/j.jestch.2019.10.006}
\end{barticle}
\endbibitem

\bibitem[\protect\citeauthoryear{Tompkins et~al.}{1996}]{tompkins1996facilities}
\begin{bbook}
\bauthor{\bsnm{Tompkins}, \binits{J.A.}},
\bauthor{\bsnm{White}, \binits{J.A.}},
\bauthor{\bsnm{Frazelle}, \binits{E.H.}},
\bauthor{\bsnm{Bozer}, \binits{Y.A.}},
\bauthor{\bsnm{Tanchoco}, \binits{J.M.A.}},
\bauthor{\bsnm{Trevino}, \binits{J.}}:
\bbtitle{Facilities Planning},
\bedition{2nd} edn.,
p. \bfpage{752}.
\bpublisher{John Wiley \& Sons Inc},
\blocation{Hoboken, NJ, USA}
(\byear{1996})
\end{bbook}
\endbibitem

\bibitem[\protect\citeauthoryear{Vizvári}{2014}]{bv2014}
\begin{botherref}
\oauthor{\bsnm{Vizvári}, \binits{B.}}:
Towards the numerical solution of a large and difficult milp problem of closed loop layout.
RUTCOR, Rutgers University, Research Report
\textbf{RRR8-2014}
(2014)
\doiurl{10.13140/2.1.1363.9366}
\end{botherref}
\endbibitem

\bibitem[\protect\citeauthoryear{Wan et~al.}{2022}]{wan2022differential}
\begin{barticle}
\bauthor{\bsnm{Wan}, \binits{X.}},
\bauthor{\bsnm{Zuo}, \binits{X.}},
\bauthor{\bsnm{Zhao}, \binits{X.}}:
\batitle{A differential evolution algorithm combined with linear programming for solving a closed loop facility layout problem}.
\bjtitle{Applied Soft Computing}
\bvolume{121},
\bfpage{108725}
(\byear{2022})
\end{barticle}
\endbibitem

\bibitem[\protect\citeauthoryear{Yang et~al.}{2023}]{YanLiuXu2023simulation}
\begin{barticle}
\bauthor{\bsnm{Yang}, \binits{C.}},
\bauthor{\bsnm{Liu}, \binits{S.}},
\bauthor{\bsnm{Xu}, \binits{Z.}}:
\batitle{A simulation-based optimization method for facility layout considering the agv path}.
\bjtitle{Journal of Physics: Conference Series}
\bvolume{2430},
\bfpage{012019}
(\byear{2023})
\doiurl{10.1088/1742-6596/2430/1/012019}
\end{barticle}
\endbibitem

\bibitem[\protect\citeauthoryear{Yang et~al.}{2005}]{Yang2005LayoutDF}
\begin{barticle}
\bauthor{\bsnm{Yang}, \binits{T.}},
\bauthor{\bsnm{Peters}, \binits{B.A.}},
\bauthor{\bsnm{Tu}, \binits{M.}}:
\batitle{Layout design for flexible manufacturing systems considering single-loop directional flow patterns}.
\bjtitle{Eur. J. Oper. Res.}
\bvolume{164},
\bfpage{440}--\blpage{455}
(\byear{2005})
\end{barticle}
\endbibitem

\bibitem[\protect\citeauthoryear{Zhang et~al.}{2017}]{zhang2017novel}
\begin{barticle}
\bauthor{\bsnm{Zhang}, \binits{B.}},
\bauthor{\bsnm{Song}, \binits{B.}},
\bauthor{\bsnm{Mao}, \binits{Z.}},
\bauthor{\bsnm{Tian}, \binits{W.}}:
\batitle{A novel wake energy reuse method to optimize the layout for savonius-type vertical axis wind turbines}.
\bjtitle{Energy}
\bvolume{121},
\bfpage{341}--\blpage{355}
(\byear{2017})
\end{barticle}
\endbibitem

\bibitem[\protect\citeauthoryear{Önüt et~al.}{2008}]{ONUT2008783}
\begin{barticle}
\bauthor{\bsnm{Önüt}, \binits{S.}},
\bauthor{\bsnm{Tuzkaya}, \binits{U.R.}},
\bauthor{\bsnm{Doğaç}, \binits{B.}}:
\batitle{A particle swarm optimization algorithm for the multiple-level warehouse layout design problem}.
\bjtitle{Computers \& Industrial Engineering}
\bvolume{54}(\bissue{4}),
\bfpage{783}--\blpage{799}
(\byear{2008})
\doiurl{10.1016/j.cie.2007.10.012}
\end{barticle}
\endbibitem

\end{thebibliography}
\end{document}